\documentclass{amsproc}

\usepackage{amsthm,amsfonts,amssymb,amsmath,amsxtra}
\usepackage[all]{xy}
\SelectTips{cm}{}
\usepackage{tikz}
\usepackage{verbatim}

\usepackage[margin=1.4in]{geometry}
\usepackage{mathrsfs}

\usepackage{mathdots}
\usepackage{MnSymbol}

\RequirePackage{xspace}
\RequirePackage{etoolbox}
\RequirePackage{varwidth}
\RequirePackage{enumitem}
\RequirePackage{tensor}
\RequirePackage{mathtools}
\RequirePackage{longtable}
\RequirePackage{multirow}
\RequirePackage{makecell}
\RequirePackage{bigints}
\RequirePackage{xcolor}

\usepackage[backref=page]{hyperref} %

\newcommand{\jiao}{\stackrel{\BL}{\cap}}
\newcommand{\Ltimes}{\stackrel{\BL}{\otimes}}

\newcommand{\sC}{\ensuremath{\mathscr{C}}\xspace}

\newcommand{\sH}{\ensuremath{\mathscr{H}}\xspace}

\newcommand{\sL}{\ensuremath{\mathscr{L}}\xspace}

\newcommand{\sP}{\ensuremath{\mathscr{P}}\xspace}

\newcommand{\BA}{\ensuremath{\mathbb{A}}\xspace}

\newcommand{\BC}{\ensuremath{\mathbb{C}}\xspace}

\newcommand{\BE}{\ensuremath{\mathbb{E}}\xspace}

\newcommand{\BG}{\ensuremath{\mathbb{G}}\xspace}

\newcommand{\BI}{\ensuremath{\mathbb{I}}\xspace}
\newcommand{\BJ}{\ensuremath{\mathbb{J}}\xspace}

\newcommand{\BL}{\ensuremath{\mathbb{L}}\xspace}

\newcommand{\BQ}{\ensuremath{\mathbb{Q}}\xspace}
\newcommand{\BR}{\ensuremath{\mathbb{R}}\xspace}
\newcommand{\BS}{\ensuremath{\mathbb{S}}\xspace}

\newcommand{\BV}{\ensuremath{\mathbb{V}}\xspace}

\newcommand{\BX}{\ensuremath{\mathbb{X}}\xspace}

\newcommand{\BZ}{\ensuremath{\mathbb{Z}}\xspace}

\newcommand{\CE}{\ensuremath{\mathcal{E}}\xspace}
\newcommand{\CF}{\ensuremath{\mathcal{F}}\xspace}

\newcommand{\CH}{\ensuremath{\mathcal{H}}\xspace}

\newcommand{\CM}{\ensuremath{\mathcal{M}}\xspace}
\newcommand{\CN}{\ensuremath{\mathcal{N}}\xspace}
\newcommand{\CO}{\ensuremath{\mathcal{O}}\xspace}

\newcommand{\CX}{\ensuremath{\mathcal{X}}\xspace}

\newcommand{\CZ}{\ensuremath{\mathcal{Z}}\xspace}

\newcommand{\RG}{\ensuremath{\mathrm{G}}\xspace}
\newcommand{\RH}{\ensuremath{\mathrm{H}}\xspace}

\newcommand{\RR}{\ensuremath{\mathrm{R}}\xspace}

\newcommand{\RU}{\ensuremath{\mathrm{U}}\xspace}

\newcommand{\RZ}{\ensuremath{\mathrm{Z}}\xspace}

 \DeclareFontFamily{U}{wncy}{}
    \DeclareFontShape{U}{wncy}{m}{n}{<->wncyr10}{}
    \DeclareSymbolFont{mcy}{U}{wncy}{m}{n}
    \DeclareMathSymbol{\Sha}{\mathord}{mcy}{"58}

\newcommand{\Ad}{\mathrm{Ad}}

\DeclareMathOperator{\Aut}{Aut}

\newcommand{\Ch}{{\mathrm{Ch}}}
\DeclareMathOperator{\charac}{char}

\newcommand{\cl}{{\mathrm{cl}}}

\newcommand{\corr}{\mathrm{corr}}

\newcommand{\del}{\operatorname{\partial Orb}}

\DeclareMathOperator{\End}{End}

\newcommand{\Mat}{\mathrm{Mat}}

\DeclareMathOperator{\Gal}{Gal}

\newcommand{\GL}{\mathrm{GL}}

\newcommand{\GU}{\mathrm{GU}}

\DeclareMathOperator{\Hom}{Hom}

\newcommand{\id}{\ensuremath{\mathrm{id}}\xspace}

\DeclareMathOperator{\Int}{\ensuremath{\mathrm{Int}}\xspace}

\DeclareMathOperator{\Lie}{Lie}

\DeclareMathOperator{\Nm}{Nm}

\DeclareMathOperator{\Orb}{Orb}
\DeclareMathOperator{\ord}{ord}

\DeclareMathOperator{\rank}{rank}
\DeclareMathOperator{\Ros}{Ros}

\newcommand{\LN}{\,^\BL\!\CN}

\DeclareMathOperator{\Pic}{Pic}

\newcommand{\rs}{\ensuremath{\mathrm{reg}}\xspace}

\newcommand{\Sh}{\mathrm{Sh}}

\newcommand{\SL}{{\mathrm{SL}}}

\DeclareMathOperator{\Spec}{Spec}
\DeclareMathOperator{\Spf}{Spf}

\newcommand{\SO}{{\mathrm{SO}}}
\newcommand{\Sp}{{\mathrm{Sp}}}

\DeclareMathOperator{\sgn}{sgn}

\DeclareMathOperator{\tr}{tr}

\newcommand{\U}{\mathrm{U}}

\newcommand{\wt}{\widetilde}
\newcommand{\wh}{\widehat}

\newcommand{\pair}[1]{\langle {#1} \rangle}

\newcommand{\ov}{\overline}

\newcommand{\incl}{\hookrightarrow}
\newcommand{\lra}{\longrightarrow}
\newcommand{\imp}{\Longrightarrow}

\newcommand{\bs}{\backslash}

\newcommand{\lv}{\lvert}
\newcommand{\rv}{\rvert}

\newlength{\olen}
\newlength{\ulen}
\newlength{\xlen}
\newcommand{\xra}[2][]{%
   \ifbool{@display}%
      {\settowidth{\olen}{$\overset{#2}{\longrightarrow}$}%
       \settowidth{\ulen}{$\underset{#1}{\longrightarrow}$}%
       \settowidth{\xlen}{$\xrightarrow[#1]{#2}$}%
       \ifdimgreater{\olen}{\xlen}%
          {\underset{#1}{\overset{#2}{\longrightarrow}}}%
          {\ifdimgreater{\ulen}{\xlen}%
             {\underset{#1}{\overset{#2}{\longrightarrow}}}
             {\xrightarrow[#1]{#2}}}}%
      {\xrightarrow[#1]{#2}}
   }

\newtheorem{theorem}{Theorem}[section]

\theoremstyle{definition}
\newtheorem{definition}[theorem]{Definition}

\theoremstyle{remark}
\newtheorem{remark}[theorem]{Remark}

\numberwithin{equation}{section}

\newtheorem{conjecture}[theorem]{Conjecture}

\newenvironment{altenumerate}
   {\begin{list}
      {(\theenumi) }
      {\usecounter{enumi}
       \setlength{\labelwidth}{0pt}
       \setlength{\labelsep}{0pt}
       \setlength{\leftmargin}{0pt}
       \setlength{\itemsep}{\the\smallskipamount}
       \renewcommand{\theenumi}{\roman{enumi}}
      }}
   {\end{list}}

   \newcommand{\isoarrow}{%
   \ifbool{@display}{\overset{\sim}{\longrightarrow}}{\xrightarrow\sim}%
   }

\newcommand{\sform}{\ensuremath{(\text{~,~})}\xspace}

\newenvironment{altitemize}
   {\begin{list}
      {$\bullet$}
      {\setlength{\labelwidth}{0pt}
	   \setlength{\itemindent}{5pt}
       \setlength{\labelsep}{5pt}
       \setlength{\leftmargin}{0pt}
       \setlength{\itemsep}{\the\smallskipamount}
      }}
   {\end{list}}

\newcommand{\iso}{\cong}
\newcommand{\mr}{\mathrm}

\newcommand{\mbE}{\mathbb{E}}
\newcommand{\mbF}{\mathbb{F}}

\newcommand{\mbV}{\mathbb{V}}

\newcommand{\mbX}{\mathbb{X}}

\newcommand{\mcE}{\mathcal{E}}

\newcommand{\mcN}{\mathcal{N}}
\newcommand{\mcO}{\mathcal{O}}

\newcommand{\mcZ}{\mathcal{Z}}

\begin{document}

\title[IHES-2022]{High dimensional Gross--Zagier formula: a survey}

\author{Wei Zhang}
\address{Massachusetts Institute of Technology, Department of Mathematics, 77 Massachusetts Avenue, Cambridge, MA 02139, USA}
\email{weizhang@mit.edu}

\thanks{Research of W.Z. is supported by NSF Grant \#1901642 and the Simons foundation.}

\date{June 10, 2022} 



\begin{abstract}We survey recent developments on generalizing the Gross--Zagier formula to high dimensional Shimura varieties, with an emphasis on the 
 Arithmetic Gan--Gross--Prasad conjecture and the relative trace formula approach.  
\end{abstract}

\maketitle

\tableofcontents
\section{Introduction}

The seminal paper  of Gross and Zagier \cite{GZ} proves a formula relating the N\'eron--Tate heights of Heegner points on modular curves to the first central derivatives of certain L-functions.
The Gross--Zagier formula and the theorem of Kolyvagin on Euler systems together imply that the (rank part of the) Birch--Swinnerton-Dyer conjecture holds if the analytic rank of a modular elliptic curve over $\BQ$ is (at most) one. In this article we survey some of the  recent developments on generalizing the Gross--Zagier formula to high dimensional Shimura varieties, as well as their applications to the generalizations of the Birch--Swinnerton-Dyer conjecture formulated by Beilinson, Bloch, and Kato.

There have been a lot of works on a parallel topic, that is, the relation between automorphic periods and special values of L-functions. Let $F$ be a global field and  $G$ a reductive group over $F$. Let $\pi$ be a (tempered) cuspidal automorphic representation of $G$. Let $H\subset G$ be a subgroup.  A fundamental question in the theory of automorphic forms is to study 
the (automorphic) $H$-period integral on $\pi$:
\begin{align}\label{eq:P pi}
\int_{[H]}\varphi(h)\, dh,\quad \varphi\in \pi.
\end{align}
Here $[H]$ denotes the automorphic quotient $H(F)\bs H(\BA)$, where $\BA$ is the ring of adeles of $F$. For many pairs $(H,G)$, the period integrals are related  
to special values of L-functions attached to $\pi$. Two notable examples are the pairs $(\GL_{n-1},\GL_{n-1}\times \GL_n)$ and $(\GL_n\times \GL_n,\GL_{2n})$. The first is related to the Rankin-Selberg convolution L-function \cite{JPSS}, the second to the standard L-function for $\GL_{2n}$ \cite{FJ} (see also \cite{BF} for a weighted version). These two examples involve {\em split} groups only; in particular,  when $F$ is a number field, the group $G$ does not satisfy the following condition (unless $G$ is commutative):
\begin{align}\label{star0}
\text{\em $G(F\otimes_\BQ\BR)$ is compact modulo its center}.
\end{align}
 The compactness of $G(F\otimes_\BQ\BR)$ can often allow us to derive algebraicity of period integrals and, together with other inputs, to relate L-values to interesting arithmetic invariants of an appropriate Galois representation associated to $\pi$. As an example in the non-split case, Waldspurger considered 
\begin{align}\label{eq:H G}
H=\RR_{F'/F}\BG_m, \quad G=B^\times,
\end{align}
where $\RR_{F'/F}$ denotes the Weil restriction of scalars, $F'$ is a quadratic field extension of $F$, with an embedding into a quaternion algebra $B$ over $F$. In the 1980s, Waldspurger \cite{W85} proved a formula of the following form (after suitably modifying the period integral to take into account the center of $G$)
$$
\left|\int_{[H]}\varphi(h)\, dh\right|^2=L(\pi_{F'}, 1/2) \times \text{``Local factors"},
$$
where $\pi_{F'}$ denotes the base change of $\pi$ to $F'$.
The group $G(F\otimes_\BQ\BR)$ is compact modulo its center if and only if $F$ is totally real and $B$ is non-split at all archimedean places. In this compact case the formula of Waldspurger played a crucial role in the study of  the  Birch--Swinnerton-Dyer conjecture for (modular) elliptic curves over a totally real number field in the analytic rank zero case (for example, in Bertolini--Darmon's variant of Kolyvagin systems \cite{BD}). 

There has been a large number of (conjectural or proved) generalizations of Waldspurger's formula to more general pairs $(H,G)$. Most relevant to this survey are the Gan--Gross--Prasad  conjecture \cite{GGP} for unitary groups and its refinement  by Ichino--Ikeda. We refer to \cite[\S2]{Z-ICM} for a survey on previous results including \cite{JR,Z14a,Z14b,Xue,BP21a} for $\pi$ satisfying certain local conditions. Recently the complete results are obtained in \cite{BPLZZ} (in the ``stable" case) and \cite{BPCZ} (in the ``endoscopic" case).   We also refer to the articles by Beuzart-Plessis \cite{BP22} and Chaudouard \cite{Ch22} on this topic. 

This survey concerns an analogous question to \eqref{eq:P pi}  where one aims to upgrade the embedding of topological spaces $[H]\subset [G]$ to one with richer structure. When $F$ is a (suitable) number field,  we expect such an upgrade when there are Shimura varieties associated to the groups $H$ and $G$. The pioneering example is the Gross--Zagier formula \cite{GZ}, proved about the same time as Waldspurger's formula. In the Gross--Zagier formula, one has, for $F=\BQ$ and $F'$ an imaginary quadratic field, 
$$
H=\RR_{F'/\BQ}\BG_m,\quad G=\GL_{2,\BQ}.
$$
Then the Gross--Zagier formula relates the N\'eron--Tate height pairing of the $\pi$-isotypical component of a special cycle (a divisor) attached to $H$ to the central derivative of the L-function for the base change $\pi_{F'}$. It applies to all $\pi$ whose local components satisfies suitable conditions. In the early 2000s, Shouwu Zhang extended the Gross--Zagier formula from modular curves to a large class of Shimura curves (over general totally real fields)  \cite{ZSW01a,ZSW01b}; see also a different generalization to Shimura curves over $\BQ$ due to Kudla--Rapoport--Yang \cite{KRY}. The most general form of the Gross--Zagier formula for Shimura curves was then obtained in the work of Xinyi Yuan, Shouwu Zhang and the author \cite{YZZ2} after another ten years. We will recall a special case in \S\ref{s:GZ};  we refer to \cite[\S3]{Z-CDM} for a survey. As mentioned already, the Gross--Zagier formula has a direct application to the  Birch--Swinnerton-Dyer conjecture. Together with Kolyvagin's method of Euler systems, it settles the rank part of the Birch--Swinnerton-Dyer conjecture for elliptic curves over $\BQ$ (more generally, modular elliptic curves over totally real fields) when the analytic rank is at most one.

After \cite{YZZ2}, a natural question is to generalize the Gross--Zagier formula for Shimura curves to higher dimensional Shimura varieties. The goal of these notes is to survey the recent developments on this topic, with an emphasis on the arithmetic Gan--Gross--Prasad (AGGP) conjecture \cite{GGP} for $\RU(n)\times \RU(n+1)$ and the relative trace formula approach initiated in \cite{Z12}.  The first proof of Waldspurger's formula mostly relies on the theory of theta series, as does the overall strategy of the proof of the Gross--Zagier formula for modular or Shimura curves in \cite{GZ,YZZ2}. Later Jacquet found a new proof of Waldspurger's formula using his relative trace formula approach, and  in \cite{JR} Jacquet and Rallis proposed a relative trace formula approach to prove the global Gan--Gross--Prasad  conjecture for unitary groups (for the automorphic period integral and central L-values). Inspired by their relative trace formula, the author formulated an approach to attack the AGGP conjecture and reduced the conjecture to several local conjectures, notably the arithmetic fundamental lemma conjecture \cite{Z09,Z12,LRZ1} and the arithmetic transfer conjectures \cite{RSZ1,RSZ2,ZZ,LRZ2}. Recent years have witnessed substantial progress in these local conjectures, see \S\ref{s:loc} for a selection of related theorems.

Here is an outline the paper. In \S\ref{s:GZ}, we recall a special case of the Gross--Zagier formula for Shimura curves over $\BQ$. 
In \S\ref{s:AGGP},  we first present the AGGP conjecture for unitary Shimura varieties and various variants (including a $p$-adic variant in \S\ref{ss:pAGGP}).  We have essentially restricted ourselves to unitary Shimura varieties with certain parahoric level structure at inert places.  One of the main difficulty is the lack of regular integral models at deeper levels. We also  recall the relative trace formula approach. In \S\ref{s:loc} we will state the local conjectures and we update the recent progress on some of the key local questions since the author's previous survey \cite{Z-ICM}: the proof of the arithmetic fundamental lemma conjecture for all $p$-adic fields (with $p$ odd) and of arithmetic  transfer conjecture in certain parahoric cases.
 In \S\ref{s:other} we will briefly mention several other generalizations of  the Gross--Zagier formula, including the AGGP conjecture in the $\RU(n)\times \RU(n)$-case formulated by Yifeng Liu, and the arithmetic Rallis inner product formula. The latter generalizes the  Gross--Zagier formula along the line of the aforementioned work of Kudla--Rapoport--Yang \cite{KRY}, and relies on the recent progress on the arithmetic Siegel--Weil formula which was envisioned by Kudla \cite{K97} and subsequently formulated as conjectures by Kudla and Rapoport in a series papers. We also present some examples of arithmetic applications to the Beilinson--Bloch--Kato conjecture in \S\ref{ss:BBK}. In another direction, there is also the higher derivative version of Gross--Zagier formula over a functional field, see \S\ref{ss:fun fd}.

 \subsection*{Acknowledgement} The author thanks Tony Feng, Yifeng Liu, Murilo Corato Zanarella for their comments on an early draft of the notes, and the anonymous referee for their careful reading and helpful comments.
 
\section{Gross--Zagier formula for Shimura curves}\label{s:GZ}
For simplicity, 
we only recall a mild generalization of the Gross--Zagier formula from modular curves $X_0(N)$ in \cite{GZ} to certain Shimura curves over $\BQ$; the ramification conditions here can all be dropped, see \cite{YZZ2}.   

Let $N>1$ be an integer. Let $K=\BQ[\sqrt{-D}]$ be an imaginary quadratic extension with fundamental discriminant $-D<0$. For simplicity we assume $D>4$ and $(D,N)=1$. This determines a unique factorization 
\begin{align}\label{eqn N+N-}
N=N^+N^-
\end{align}
 where the prime factors of $N^+$ ($N^-$, resp.) are all split (inert, resp.) in $K$. Let us impose the {\em generalized Heegner hypothesis}: 
 $$
 \text{ $N^-$ is square-free with an {\em even} number of prime factors.}
 $$

 We consider the indefinite quaternion algebra $B$ over $\BQ$ that is ramified precisely at all factors of $N^-$. Let $G=B^\times$ viewed as an algebraic group over $\BQ$. We have the usual Shimura datum attached to $G$. Consider the (compactified)  Shimura curve $X_{N^+,N^-}=X_U$ where the compact open $U\subset G(\BA_f^\times)=B^\times(\BA_f)$ is defined by
$$U_{N^+,N^-}=\prod_{\ell<\infty}U_\ell,\quad U_\ell=\begin{cases}\Gamma_0(N),& \ell|N^+,\\
\CO_{B_\ell}^\times,& \ell\nmid N^+.
\end{cases}$$
Equivalently, we may consider an Eichler order $\CO_{B,N^+}$ in $\CO_B$ with level $N^+$ and define 
\begin{align}\label{eqn shi curve}
U_{N^+,N^-}=(\CO_{B,N^+}\otimes\wh{\BZ})^\times.
\end{align}
We have an isomorphism for the complex analytic space
\begin{align}\label{eqn X C}
X_{N^+,N^-}(\BC)=B(\BQ)\bs (\CH \times B(\BA_f))/(\CO_{B,N^+}\otimes\wh{\BZ})^\times \cup\{\text{cusps}\}.
\end{align}
In particular, if $N^-=1$, the curve $X_{N^+,N^-}$ is the  (compactified) classical modular curve $X_0(N^+)$, whose complex points are $\Gamma_0(N^+)\bs \CH$ together with cusps.

The case $N^-=1$ corresponds to the  classical {\em Heegner hypothesis} in \cite{GZ}: {\em every prime factor $\ell|N$ is split in $K$.} Then there 
exists an ideal $\CN$ of $\CO_K$, the ring of integers of $K$, such that
$$
\CO_K/\CN\simeq\BZ/N\BZ.
$$
Then the elliptic curves $\BC/\CO_K$ and $\BC/\CN^{-1}$ are naturally related by a cyclic isogeny of degree $N$, therefore define a $\BC$-point, denoted by $x(1)$, on $X_0(N)$. By the theory of complex multiplication, the point $x(1)$ is defined over the Hilbert class field $H$ of $K$. (This depends on a (non-unique) choice of an ideal $\CN$ of $\CO_K$ above; but they all yield the same Gross--Zagier formula below.) By class field theory, the Galois group $\Gal(H/K)$ is isomorphic to the class group $\Pic(\CO_K)$. Similarly, under the generalized Heegner hypothesis, we can define a point $x(1)\in X_{N^+,N^-}(H)$. 

Let $E/\BQ$ be an elliptic curve with conductor $N$. Consider a modular parameterization 
 $$
 \phi: X_{N^+,N^-}\lra E.
$$If $N^-=1$ we require this morphism to send $(\infty)$ to $0$ of $E$; in general we require that it pushes forward the class of the Hodge bundle (in the Picard group) to a multiple of the class of $0$ of $E$.
Define
\begin{align}\label{eqn def y1 yK}
y(1)=\phi(x(1))\in E(H),\quad y_K=\tr_{H/K}\,y(1)\in E(K).
\end{align}

Let $f$ be the newform of weight two and $\Gamma_0(N)$-level associated to the elliptic curve $E$.
Let $L(f/K,s)=L(E/K,s)$ be the L-function (without the archimedean factor) attached to the elliptic curve $E$ base changed to $K$; here we take the classical normalization for $L(f/K,s)$, i.e, the center of functional equation is at $s=1$. The generalized Heegner hypothesis ensures that the global root number $$\epsilon(E/K)=-1.$$ Then the central value $L(f/K,1)=0$ and the Gross--Zagier formula for $X_{N^+,N^-}$ is as follows:

\begin{align}\label{eqn GZ f/K}
\frac{L'(f/K,1)}{(f,f)}=\frac{1}{\sqrt{|D|}}\frac{\pair{y_K,y_K}_{\rm NT}}{\deg(\phi)},
\end{align}
where $\pair{y_K,y_K}_{\rm NT}$ is the N\'eron--Tate height pairing over $K$ and $(f,f)$ the Petersson inner product \begin{align}\label{eqn def (f,f)}
(f,f):=8\pi^2\int_{\Gamma_0(N)\bs\CH} f(z)\ov{f(z)}dxdy=\int_{X_0(N)(\BC)}\omega_f\wedge \ov{i\omega}_f,
\end{align}
where $\omega_f:=2\pi i f(z)dz$.

This is proved by Gross--Zagier in \cite{GZ} when $N^-=1$, in general by S. Zhang in \cite{ZSW01a,ZSW01b} and Yuan--Zhang--Zhang in \cite{YZZ2}.

\begin{remark}To compare with
 the Birch--Swinnerton-Dyer conjectural formula for $L'(E/K,1)$, we consider a modular parameterization, still denoted by $f$,
\begin{align}\label{eqn mod para}
f: X_0(N)\to E,
\end{align}which we assume maps the cusp $(\infty)$ to zero.
The pull-back of the N\'eron differential $\omega$ on $E$ is $f^\ast \omega=c\cdot \omega_f$ for a constant $c$. It is known that the constant $c$ is an integer. If $f$ is an optimal  parameterization, it is called the Manin constant and conjecturally equal to $1$. By a theorem of Mazur, we have $p|c\imp p|2N$ (there are various improvements of Mazur's theorem but for the sake of brevity we do not recall them).
Then an equivalent form of the formula (\ref{eqn GZ f/K}) is
\begin{align}\label{eqn GZ E/K}
\frac{L'(E/K,1)}{|D|^{-1/2} \int_{E(\BC)}\omega\wedge \ov{\omega}}=\frac{1}{c^2}\frac{\deg(f)}{\deg(\phi)}\pair{y_K,y_K}_{\rm NT}.
\end{align}
For comparison, the Birch--Swinnerton-Dyer conjecture asserts that 
\begin{align}\label{eqn BSD E/K}
\frac{L'(E/K,1)}{|D|^{-1/2} \int_{E(\BC)}\omega\wedge \ov{\omega}}=\frac{ \prod_{v|N}m_v}{[E(K): \BZ y_K]^2}\cdot\# \Sha_{E/K} \cdot \pair{y_K,y_K}_{\rm NT},
\end{align}
where $m_v$ denotes the local Tamagawa number at the bad places of $E/K$, $[E(K): \BZ y_K]$ denotes the index of the subgroup generated by $y_K\in E(K)$, and $\# \Sha_{E/K} $ denotes the order of the Tate--Shafarevich group.  Therefore, by \eqref{eqn GZ E/K} the conjectural formula \eqref{eqn BSD E/K} is equivalent to 
$$
\# \Sha_{E/K}=\frac{[E(K): \BZ y_K]^2 \cdot (\deg(f)/\deg(\phi))}{c^2 \cdot \prod_{v|N}m_v}.
$$
\end{remark}

\section{AGGP conjecture for unitary Shimura varieties}\label{s:AGGP}
We note that from now on we may use different notation from the previous section.
 \subsection{Special pairs of Shimura data and special cycles}
\label{sss:pair Sh}
This subsection is largely from \cite{Z-ICM}. We recall the notion of {\em a special pair of Shimura data},  which is an enhancement of a spherical pair $(\RH, \RG)$ to a homomorphism of Shimura data $(\RH,X_\RH)\to (\RG,X_\RG)$. 

Let  $\BS$ be the torus $\RR_{\BC/\BR}\BG_m$ over $\BR$ (i.e., we view $\BC^\times$ as an $\BR$-group). Recall that a Shimura datum $\bigl(\RG,X_\RG\bigr)$ consists of a reductive group $\RG$ over $\BQ$, and a $\RG(\BR)$-conjugacy class $X_\RG=\{h_\RG\}$ of $\BR$-group homomorphisms $h_\RG:\BS\to \RG_\BR$ (sometimes called Shimura homomorphisms)  satisfying Deligne's list of axioms \cite[1.5]{De}. In particular, $X_\RG$ is a Hermitian symmetric domain. 
 \begin{definition}
 A {\it special pair} of Shimura data is a homomorphism  \cite[1.14]{De} between two Shimura data 
 $$
 \xymatrix{\delta\colon\bigl(\RH,X_\RH\bigr)\ar[r]& \bigl(\RG,X_\RG\bigr) }
 $$ such that
  \begin{altenumerate}
 \item the homomorphism $\delta: \RH\to \RG$ is injective such that the pair $(\RH,\RG)$ is spherical, and 
 \item the dimensions of $X_\RH$ and $X_\RG$ (as complex manifolds) satisfy 
 $$
\dim_\BC X_\RH=\bigg\lfloor \frac{\dim_\BC X_\RG}{2}\bigg\rfloor.
 $$ 
 \end{altenumerate}
 \end{definition} 
  \begin{remark}
The condition on the pair $(\RH,\RG)$ being spherical  may be an artifact. There are certainly examples of non-spherical pairs that result to Shimura data satisfying the two conditions above. However, all the ``interesting" (as of today) examples are spherical and there are already many of them. Moreover,  interesting automorphic period integrals seem mostly (if not all) attached to spherical pairs.
\end{remark}

 \begin{remark}
 It seems an interesting question to enumerate special pairs of Shimura data.  In fact, we may consider the analog of special pairs of Shimura data in the context of {\em local Shimura data}. It seems more realistic to enumerate the pairs in the local situation. Moreover, it is only a matter of interest in this paper that we impose the dimension equality in (ii). There are certainly other interesting cases where the equality does not hold. 

\end{remark}

For a Shimura datum $  (\RG,X_\RG)$ we have a projective system of {\em Shimura varieties} $\{\Sh_{K}(\RG)\}$, indexed by compact open subgroups $K\subset\RG(\BA_f)$,  of smooth quasi--projective varieties (for neat $K$) defined over a number field $E$---the reflex field of $(\RG,X_\RG)$, viewed as a subfield of $\BC$.

For a special pair of Shimura data $(\RH,X_\RH)\to (\RG,X_\RG)$, compact open subgroups $K_\RH\subset \RH(\BA_f)$ and $K_\RG\subset\RG(\BA_f)$ such that $K_\RH\subset K_\RG$, we have a finite morphism (over the reflex field of $(\RH,X_\RH)$)
$$
\delta_{K_\RH,K_\RG}:\Sh_{K_\RH}(\RH)\to \Sh_{K_\RG}(\RG).
$$
The cycle $z_{K_\RH,K_\RG}:=\delta_{K_\RH,K_\RG,\ast} [\Sh_{K_\RH}(\RH)]$ on $ \Sh_{K_\RG}(\RG)$ will be called the {\em special cycle} (for the level $(K_\RH,K_\RG)$).  Very often we choose $K_\RH=K_\RG\cap \RH(\BA_f)$ in which case we simply denote the special cycle by $z_{K_\RG}$.

There are two cases depending on the parity of the dimension of $ \Sh_{K_\RG}(\RG)$:
\begin{itemize}
\item
When $\dim_\BC X_\RG$ is even, the special cycles are in the middle dimension.  They can very often account for most Tate cycles on $ \Sh_{K_\RG}(\RG)$ and are therefore related to poles of L-functions at the near central point $s=1$. They are related to automorphic period integrals and distinguished automorphic representations as in many works of Jacquet. One example is due to Harder, Langlands, and Rapoport \cite{HLR} on Tate cycles on Hilbert modular surfaces, which has led to  many  generalizations.
\item When $\dim_\BC X_\RG$ is odd, the special cycles are just below the middle dimension, and we will say that they are in the {\em arithmetic  middle dimension} (in the sense that, once extending both Shimura varieties to suitable integral models, we obtain cycles in the middle dimension). Their height pairings are expected to be related to  the central derivative of L-functions and therefore can be used to attack  Beilinson and Bloch's generalizations of  Birch--Swinnerton-Dyer conjecture to general varieties, at least when the relevant L-function has vanishing order one. 

\end{itemize}

Below we focus on the case where the special cycles are in the arithmetic  middle dimension. 

\subsection{The main examples for this survey}

We will focus on the Gan--Gross--Prasad pair for unitary groups; for more examples see \S\ref{s:other}.
 \subsubsection{Gross--Zagier pair}
  In the case of the Gross--Zagier formula \cite{GZ}, one considers an embedding of an {\em imaginary} quadratic field $F'$ into $\Mat_{2,\BQ}$ (the $\BQ$-algebra  of $2\times 2$-matrices), and the induced embedding  
  $$\RH=\RR_{F'/\BQ}\BG_m\incl \RG=\GL_{2,\BQ}.$$   Note that $\RH_\BR\simeq\BC^\times$ as $\BR$-groups (upon a choice of embedding $F'\incl \BC$). This defines  $h_\RH: \BS\to \RH_\BR$, and its composition with the embedding $ \RH_\BR\to \RG_\BR$ defines $h_{\RG}:\BS\to \BG_\BR$. We obtain a  special pair $(\RH,X_\RH)\to (\RG,X_\RG)$, 
where
  $$
  \dim  X_\RG=1,\quad\quad \dim  X_\RH=0.
  $$
  
  In the general case, we replace $F'/\BQ$ by a CM extension $F'/F$ of a totally real number field $F$, and replace $\Mat_{2,\BQ}$ by a quaternion algebra $B$ over $F$ that is ramified at all but one archimedean places of $F$. In this generality the Gross--Zagier formula is treated in \cite{YZZ2}.

 \subsubsection{Gan--Gross--Prasad pair}
  \label{sss:GGP pair}
Let $F$ be a number field, and let $F'=F$ in the orthogonal case and $F'$ a quadratic extension of $F$ in the Hermitian case. Let $W_n$ be a non-degenerate orthogonal space or Hermitian space with $F'$-dimension $n$.  Let $W_{n-1}\subset W_n$ be a non-degenerate subspace of codimension one. Let $\RG_{i}$ be $\SO(W_i)$ or $\RU(W_i)$ for $i=n-1,n$, and $\delta:\RG_{n-1}\incl\RG_n$ the induced embedding. Let 
\begin{align}
 \RG=\RG_{n-1}\times \RG_{n},
	\quad\,\quad
	\RH=\RG_{n-1},
\end{align} with the ``diagonal" embedding $\Delta:\RH \incl \RG$ (i.e., the graph of $\delta$). The pair $(\RH,\RG)$ is spherical and we call it the {\em Gan--Gross--Prasad pair}.

 Now we impose the following conditions.
 \begin{altenumerate}
 \item $F$ is a totally real number field, and in the Hermitian case $F'/F$ is a CM (=totally imaginary quadratic) extension.
 \item  For an archimedean place $\varphi\in\Hom(F,\BR)$, denote by $\sgn_{\varphi}(W)$ the signature of $W\otimes_{F,\varphi}\BR$ as an orthogonal or Hermitian space over $F'\otimes_{F,\varphi}\BR$. Then there exists a  distinguished real place $\varphi_0\in\Hom(F,\BR)$ such that
 $$
 \sgn_{\varphi}(W_{n})=\begin{cases}(2,n-2),&\varphi=\varphi_0
 \\(0,n),& \varphi\in\Hom(F,\BR)\setminus\{\varphi_0\}
 \end{cases}
 $$
 in the orthogonal case, and 
  $$
 \sgn_{\varphi}(W_{n})=\begin{cases}(1,n-1),&\varphi=\varphi_0
 \\(0,n),& \varphi\in\Hom(F,\BR)\setminus\{\varphi_0\}
 \end{cases}
 $$in the Hermitian case. In addition, the quotient $W_n/W_{n-1}$ is negative definite at every $\varphi\in\Hom(F,\BR)$ (so the signature of $W_{n-1}$ is given by similar formulas).
\end{altenumerate} Then Gan, Gross, and Prasad \cite[\S27]{GGP} define Shimura data that enhance
the embedding $\RH\incl \RG$ to a homomorphism of Shimura data $(\RH,X_\RH)\to (\RG,X_\RG)$, where the dimensions are 
$$\begin{cases}
\dim  X_\RG=2n-5,\quad \dim  X_\RH=n-3,\quad \text{in the orthogonal case,}\\
\dim  X_\RG=2n-3,\quad \dim  X_\RH=n-2,\quad \text{in the Hermitian case}.
\end{cases}
$$ 
The reflex field is $F'$  ($=F$ in the orthogonal case) via the distinguished embedding $\varphi_0$.

\subsection{The arithmetic Gan--Gross--Prasad conjecture}

\subsubsection{Height pairings}
\label{sss std conj}

 Let $X$ be a smooth proper variety over a number field $E$, and let $\Ch^i(X)$ be the group of codimension-$i$ algebraic cycles on $X$ modulo rational equivalence. We have a cycle class map
$$
   \cl_i\colon \Ch^i(X)_\BQ\to H^{2i}(X),
$$where $H^{2i}(X)$ is the Betti cohomology
$   H^*\bigl(X(\BC),\BC\bigr).$
The kernel is the group of cohomologically trivial cycles, denoted by $\Ch^i(X)_{0}$.

 Conditional on certain standard conjectures on algebraic cycles, there is a height pairing defined by Beilinson and Bloch, 
 \begin{equation}
\label{eqn BB}
  \sform_\mathrm{BB}\colon \Ch^i(X)_{\BQ,0}\times \Ch^{d+1-i}(X)_{\BQ,0}\to \BR,\quad d=\dim X.
\end{equation}

This is unconditionally defined when $i=1$ (the N\'eron--Tate height), or when $X$ is an abelian variety \cite{Kun2}. In general, there are several approaches to the definitions (conjecturally they all yield the same height pairing). In one of the approaches, if there exists a regular integral model and the cycles have an appropriate extension, so-called ``flat extension", or ``admissible extension" (\`a la Shouwu Zhang) to the integral models, then the height pairing can be defined  in terms of the arithmetic intersection theory of Arakelov and Gillet--Soul\'e \cite[\S4.2.10]{GS}, see the paper of K\"unnemann \cite{Kun1} for a precise statement. In particular, when there exists a smooth proper model $\CX$ of $X$ over $O_E$ (this is also true for Deligne--Mumford (DM) stacks $X$ and $\CX$), cf.\ \cite[\S6.1]{RSZ3}, one can define the height pairing unconditionally. In general it remains a tantalizing question how to make an unconditional height pairing.

\subsubsection{The arithmetic Gan--Gross--Prasad conjecture}\label{ss: AGGP}

We consider the special cycle in the Gan--Gross--Prasad setting \S\ref{sss:GGP pair}, which we also call the {\em arithmetic diagonal cycle} \cite{RSZ3}. We will state a version of the  arithmetic Gan--Gross--Prasad conjecture assuming some standard conjectures on algebraic cycles (cf.\ \cite[\S6]{RSZ3}), in particular, that we have the height pairing \eqref{eqn BB}. 

For each $K\subset \RG(\BA_f)$, one can construct ``Hecke--Kunneth" projectors that project the total cohomology of the Shimura variety $\Sh_{K}(\RG)$ (or its toroidal compactification) to the odd-degree part (cf.\ \cite[\S6.2]{RSZ3} in the Hermitian case; the same proof works in the orthogonal case).  Then we apply this projector to define a cohomologically trivial cycle $z_{K,0}\in \Ch^{n-1}\bigl(\Sh_{K}(\RG)\bigr)_{0}$ (with $\BC$-coefficient). Note that if we assume certain conjectures of Beilinson on algebraic cycles over number fields (cf.\ \cite[Remark 6.9]{RSZ3}), the classes $\{z_{K,0}\}_{K\subset \RG(\BA_f)}$ are independent of the choice of our projectors and they form a projective system (with respect to push-forward).  

We form the colimit $$\Ch^{n-1}\bigl(\Sh(\RG)\bigr)_{0}:=\varinjlim_{K\subset \RG(\BA_f)} \Ch^{n-1}\bigl(\Sh_{K}(\RG)\bigr)_{0}.$$ The (conditionally defined) height pairing \eqref{eqn BB}  with $\{z_{K,0}\}_{K\subset \RG(\BA_f)}$ defines a linear functional
$$
\xymatrix{\sP_{\Sh(\RH)}\colon \Ch^{n-1}\bigl(\Sh(\RG)\bigr)_{0}\ar[r]&\BC}.
$$
  This is the arithmetic version of the automorphic period integral \eqref{star0}. The group $ \RG(\BA_f)$ acts on the space $ \Ch^{n-1}\bigl(\Sh(\RG)\bigr)_{0}$. For any representation $\pi_f$ of $ \RG(\BA_f)$, let $\Ch^{n-1}\bigl(\Sh(\RG)\bigr)_{0}[\pi_f]$ denote the $\pi_f$-isotypic component of the Chow group $\Ch^{n-1}\bigl(\Sh(\RG)\bigr)_{0}$. Conjecturally the subspace $\Ch^{n-1}\bigl(\Sh(\RG)\bigr)_{0}[\pi_f]$ should be an admissible representation of $ \RG(\BA_f)$ and $\Ch^{n-1}\bigl(\Sh(\RG)\bigr)_{0}$ should be a direct sum of  $\Ch^{n-1}\bigl(\Sh(\RG)\bigr)_{0}[\pi_f]$  for $\pi_f$ appearing in the middle degree cohomology $H^{2n-3}\bigl(\Sh(\RG),\BC\bigr)$; see \cite[\S6, Remark 6.9]{RSZ3}.
  
  We are ready to state the arithmetic Gan--Gross--Prasad conjecture \cite[\S27]{GGP}, parallel to  the global GGP conjecture for the central L-value.
\begin{conjecture}\label{conj AGGP}
Let $\pi$ be a tempered cuspidal automorphic
representation of $\RG(\BA)$, appearing in the cohomology $H^{\ast}(\Sh(\RG))$.
 The following statements are equivalent.
\begin{altenumerate}
\item The linear functional $\sP_{\Sh(\RH)}$ does not vanish on the $\pi_f$-isotypic component $\Ch^{n-1}\bigl(\Sh(\RG)\bigr)_{0}[\pi_f]$.
\item  The space $\Hom_{\RH(\BA_f)}(\pi_f,\BC)\neq 0$ and the first order derivative  $L'(\frac{1}{2},\pi,R)\neq 0$. 
\end{altenumerate}
\end{conjecture}
Here $L(s,\pi,R)$ is the Rankin--Selberg convolution L-function of the base change of the two factors of $\pi=\pi_{n-1}\boxtimes\pi_n$.

 In the orthogonal case with $n\leq 4$, and when the ambient Shimura variety is a curve ($n=3$), or a product  of three curves ($n=4$),  the conjecture is unconditionally formulated. The case $n=3$  is proved by X. Yuan, S. Zhang, and the author in \cite{YZZ2}; in fact we proved a refined version.  When $n=4$ and in the triple product case (i.e., the Shimura variety $\Sh_{K}(\RG)$ is a product of three curves), X. Yuan, S. Zhang, and the author formulated a refined version of the above conjecture and proved it in some special cases, cf.\ \cite{YZZ3}.

Since the structure of $\Ch^{n-1}\bigl(\Sh(\RG)\bigr)_{0}$ as a representation of $ \RG(\BA_f)$ depends on widely open conjectures on algebraic cycles over number fields, there seems no hope to prove the above conjecture. Therefore it is desirable to formulate more accessible ones. 
Fix $K=\prod_{v}K_v\subset \RG(\BA_f)$ and consider the Hecke algebra $\sH( \RG,K)$ of bi-K-invariant functions (valued in $\BC$). Let $S$ be a large finite set such that $K_v$ is hyperspecial for all $v\notin S$ and let $\sH( \RG,K^S)=\otimes_{v\notin S}\sH( \RG,K_v) $ denote the away-from-$S$ spherical Hecke algebra. Let $E$ denote the reflex field and let $R(f)\in \Ch^{2n-3}(\Sh_K(\RG)\times_{\Spec E} \Sh_K(\RG))$ denote the associated Hecke correspondence. Define 
  \begin{align}\label{eqn BB K}
\Int_{\rm BB}(f):=(R(f)\ast z_{K,0}, z_{K,0})_\mathrm{BB}
\end{align}
whenever the right hand side is defined (e.g., by verifying the conditions in \cite{Kun1}).

\begin{conjecture}\label{conj:sph AGGP}[Assume that the height pairing \eqref{eqn BB K} is well-defined]
There is a decomposition $$
 \Int_{{\rm BB}}(f)=\sum_{\pi}\Int_\pi(f),   \quad\text{for all }\, f\in \sH( \RG,K),
 $$such that 
 \begin{enumerate}
 \item
  the sum runs over all cohomological (with respect to the trivial coefficient system) automorphic representations of $\RG(\BA)$ with non-trivial $K$-invariants, that appear in the middle-degree cohomology $H^{2n-3}(\Sh_{ K}(\RG))$;
  \item  $\Int_\pi$ is an eigen-distribution for the spherical Hecke algebra $\sH(\RG,K^S) $  with eigen-character $\lambda_{\pi}$ 
  in the sense that for all $f_0 \in \sH\left(\RG, K\right)$ and $f\in \sH(\RG,K^S)$,
\begin{align}\label{eigen dis}
\Int_\pi(f\ast f_0)=\lambda_{\pi}\left(f\right)\,\Int_\pi\left(f_0\right),
\end{align}
where $\lambda_{\pi}$ is the ``eigen-character" of $ \sH\left(\RG, K\right)$ associated to $\pi$.
\end{enumerate}
 
 Moreover, if such a representation $\pi$ is tempered, then 
 \begin{align}\label{eq:AGGP ref}
 \Int_\pi(f)=2^{-\beta_\pi} \sL'(1/2,\pi) \prod_{v<\infty}\BI_{\pi_v}(f_v).
 \end{align}
\end{conjecture}

We explain the undefined terms above. 
Here $ \sL(s,\pi)$ is the L-functions appearing in the GGP conjecture:
\begin{align}
\sL(s,\pi)=\Delta_{n}\frac{L(s,\pi,R)}{L(1,\pi, \Ad)},
\end{align}
where $\Delta_{n}$ is the special value of a certain Artin L-function (attached to the ``Gross motive" of $G_n$). 
We also write $\sL(s,\pi_v)$ for the corresponding local factor at a place $v$. The constant $\beta_\pi$ is a power of $2$ depending on the number of cuspidal factors in the isobaric sum of the base change of $\pi$.  The local spherical characters $\BI_{\pi_v}$ for the triple $( \RG,\RH, \RH)$  are defined as
 \begin{align*}
\BI_{\pi_v}(f_v):=\sum_{\phi_v\in {\rm OB}(\pi_v)} \alpha_v(\pi_v(f_v)\phi_v,\phi_v),\quad f_v\in \sC^\infty_c(\RG(F_v)),
\end{align*}
where the sum runs over an orthonormal basis $ {\rm OB}(\pi_v)$ of $\pi_v$ and where $ \alpha_v$ is the normalized local canonical invariant form of Ichino--Ikeda \cite{II}, an integration of matrix coefficient of the unitary tempered representation $\pi_v$ with invariant inner product ${\langle-,  -
\rangle _{\pi_v}} $:
\begin{align}\label{alpha1}
\alpha_v(\phi_v,\varphi_v)=\frac{1}{\sL(\frac{1}{2},\pi_v)}\int_{\RH(F_v)}{\langle \pi_v(h)\phi_v,\varphi_v
\rangle _{\pi_v} }\, dh.
\end{align}

 The formula \eqref{eq:AGGP ref} is a generalization of the Gross--Zagier formula, and should be compared with the refinement due to Ichino--Ikeda of the global Gan--Gross--Prasad conjecture, reformulated in terms of  local spherical characters (see \cite[Conj. 1.6]{Z14b}).

 Now the hope is to find as many $K$ and functions $f\in\sH( \RG,K)$
as we can such that the height pairing \eqref{eqn BB K} is defined. Along this line, we formulate unconditional conjectures below for some special levels $K \subset\RG(\BA_f)$ for unitary groups.

\subsection{The case of smooth integral models}

\label{sss: sph AGGP}

In fact in  \cite{RSZ3} Rapoport, Smithling, and the author work with a variant of the GGP Shimura data.  We modify the groups $\RR_{F/\BQ}\RG$ and $\RR_{F/\BQ}\RH$ defined previously
{\allowdisplaybreaks
\begin{align*}\RZ^\BQ&:=\GU_1=\bigl\{\,   z\in  \RR_{F'/\BQ}\BG_m \bigm| \Nm_{F'/F}(z) \in \BG_m  \,\bigr\} ,\\
\wt \RH &:={\rm G}\left(\U_1\times \U(W_{n-1})\right)=\bigl\{\, (z,h) \in  \RZ^\BQ\times \GU(W_{n-1})\bigm| \Nm_{F'/F}(z) = c(h) \,\bigr\},\\
\wt \RG &:={\rm G}\left(\U_1\times \U(W_{n-1})\times \U(W_{n}) \right)\\
&\phantom{:}= \bigl\{\, (z,h,g) \in \RZ^\BQ \times \GU(W_{n-1})\times  \GU(W_{n})\bigm| \Nm_{F'/F}(z) = c(h)=c(g) \,\bigr\},
\end{align*}
where the symbol $c$ denotes the unitary similitude factor. 
Then we have 
\begin{equation}\label{proddec}
   \begin{gathered}
   \begin{gathered}
	\xymatrix@R=0ex{
      \wt \RH \ar[r]^-\sim  &  \RZ^\BQ \times \RR_{F/\BQ} \RH
	}
   \end{gathered},
	\qquad
   \begin{gathered}
	\xymatrix@R=0ex{
	   \wt \RG \ar[r]^-\sim  & \RZ^\BQ \times \RR_{F/\BQ} \RG	}
	\end{gathered}.
   \end{gathered}
\end{equation}

We then define natural Shimura data $\big(\wt \RH,\{h_{\wt \RH}\} \big)$ and $\big(\wt \RG,\{h_{\wt \RG}\}\big)$, cf.\ \cite[\S3]{RSZ3}. This variant has the nice feature that the Shimira varieties are of PEL type, i.e., the canonical models are related to moduli problems of abelian varieties with polarizations, endomorphisms, and level structures, cf.\ \cite[\S4--\S5]{RSZ3}. 

For suitable Hermitian spaces and a special level structure $ K_{\wt \RG}^\circ\subset \wt\RG(\BA_f)$, we can even define {\em smooth} integral models (over the ring of integers of the reflex field) of  the Shimura variety $\Sh_{ K_{\wt \RG}^\circ}(\wt \RG)$. For a general CM extension $F'/F$, it is rather involved to state this level structure \cite[Remark 6.19]{RSZ3} and define the integral models \cite[\S5]{RSZ3}. For simplicity, from now on we consider a special case, when $F=\BQ$ and $F'=F[\varpi]$ is an imaginary quadratic field. We further assume that the prime $2$ is split in $F'$. We choose $\varpi\in F'$ such that $(\varpi)\subset O_{F'}$ is the product of all ramified prime ideals in $O_{F'}$.  

We first define an auxiliary moduli functor $\CM_{(r,s)}$ over $\Spec O_{F'}$ for $r+s=n$ (similar to \cite[\S13.1]{KR-U2}). 
For a locally noetherian scheme $S$ over $\Spec O_{F'}$, $\CM_{(r,s)}(S)$ is the groupoid of triples
 $(A,\iota,\lambda)$
where  
 \begin{altitemize}
\item[(i)] $(A,\iota)$ is an abelian scheme over $S$, with $O_{F'}$-action $\iota: O_{F'}\to \End(A)$ satisfying the Kottwitz condition of signature $(r, s)$, and 
\item[(ii)] $\lambda:A\to A^\vee$ is a polarization whose Rosati involution induces on $O_{F'}$ the non-trivial Galois automorphism of $F'/F$, and such that $\ker(\lambda)$ is contained in $ A[\iota(\varpi)]$ of rank $\#(O_{F'}/(\varpi))^{n}$ (resp.\ $\#(O_{F'}/(\varpi))^{n-1}$) when $n=r+s$ is even (resp.\ odd).
In particular, we have $\ker(\lambda)=A[\iota(\varpi)]$ if $n=r+s$ is even.
\end{altitemize}
 Now we assume that $(r,s)=(1,n-1)$ or $(n-1,1)$. We further impose the {\em wedge condition} and the {\em (refined) spin condition}, cf.\ \cite[\S4.4]{RSZ3}. The functor is represented by a Deligne--Mumford stack again denoted by $\CM_{(r,s)}$.
It is {\em smooth} over $\Spec O_{F'}$, despite the ramification of the field extension $F'/\BQ$, cf.\ \cite[\S4.4]{RSZ3}. 
Then we have an integral model of copies of the Shimura variety $\Sh_{ K_{\wt \RG}^\circ}(\wt \RG)$ defined by
$$
\CM_{ K_{\wt \RG}^\circ}(\wt \RG)=\CM_{(0,1)} \times_{\Spec O_{F'}} \CM_{(1,n-2)}\times_{\Spec O_{F'}}\CM_{(1,n-1)}.
$$  (In \cite[\S5.1]{RSZ3} we do cut out the desired Shimura variety with the help of a sign invariant. Moreover, since we are working with $F=\BQ$ here, we implicitly need to replace this space by its toroidal compactification.)

We now describe the arithmetic diagonal cycle (or rather, its integral model) for the level $K_{\wt \RH}^\circ=K_{\wt \RG}^\circ\cap \wt \RH(\BA_f)$. When $n$ is odd (so $n-1$ is even), we define 
$$
 \CM_{K_{\wt \RH}^\circ} \bigl(\wt \RH\bigr)=\CM_{(0,1)} \times_{\Spec O_{F'}} \CM_{(1,n-2)},
 $$ and we can define an embedding explicitly by ``taking products" (one sees easily that the conditions on the kernels of polarizations are satisfied):
\begin{equation}\label{inc M}
   \xymatrix@R=0ex{
	   \CM_{K_{\wt \RH}^\circ} \bigl(\wt \RH\bigr) \ar[r]  & \CM_{K_{\wt\RG}^\circ}\bigl(\wt \RG\bigr)\\
		\bigl( A_0,\iota_0,\lambda_0,A^\flat,\iota^\flat,\lambda^\flat \bigr) \ar@{|->}[r]
		   & \bigl( A_0,\iota_0,\lambda_0,A^\flat,\iota^\flat,\lambda^\flat, A^\flat \times  A_0,\iota^\flat \times \iota_0,\lambda^\flat \times \lambda_0\bigr) .
	}
\end{equation}	
When  $n$ is even, the situation is more subtle; see \cite[\S4.4]{RSZ3}.

 With the smooth integral model, we have an unconditionally defined height pairing \eqref{eqn BB} on $X=\Sh_{ K_{\wt \RG}^\circ}(\wt \RG)$. Now we again apply a suitable Hecke--Kunneth projector to the cycle $z_K$ for $K=K_{\wt \RG}^\circ$, and we obtain a cohomologically trivial cycle $z_{K,0}\in\Ch(\Sh_{ K_{\wt \RG}^\circ}(\wt \RG))_0 $. Then the pairing \eqref{eqn BB K} is well-defined
  \begin{align}\label{int f}
 \Int_{{\rm BB}}(f)=\Big( R(f)\ast z_{K,0},\,\,z_{K,0}\Big)_{{\rm BB}},\quad f\in \sH\left(\wt\RG, K_{\wt \RG}^\circ\right).
 \end{align}

We can restate Conjecture \ref{conj:sph AGGP} as follows.

\begin{conjecture}\label{conj:sph AGGP sm}
There is a decomposition $$
 \Int_{{\rm BB}}(f)=\sum_{\pi}\Int_\pi(f), \quad\text{for all }\, f\in \sH\left(\wt\RG, K_{\wt \RG}^\circ\right),
 $$such that 
 \begin{enumerate}
 \item
  the sum runs over all cohomological (with respect to the trivial coefficient system) automorphic representations of $\wt\RG(\BA)$ with non-trivial $K_{\wt \RG}^\circ$-invariants, that appear in the middle-degree cohomology $H^{2n-3}(\Sh_{ K_{\wt \RG}^\circ}(\wt \RG))$ and are  trivial on $Z^\BQ(\BA)$ (hence such $\pi$ essentially  comes from an automorphic representation of $U(W_{n-1})\times U(W_{n})$);
  \item  $\Int_\pi$ is an eigen-distribution for the spherical Hecke algebra $\sH(\wt \RG,K_{\wt \RG}^\circ) $  with eigen-character $\lambda_{\pi}$ 
  in the sense that for all $f \in \sH\left(\wt\RG, K_{\wt \RG}^\circ\right)$,
\begin{align}\label{eigen dis}
\Int_\pi(f)=\lambda_{\pi}\left(f\right)\,\Int_\pi\left(1_{K_{\wt \RG}^\circ }\right),
\end{align}
where $\lambda_{\pi}$ is the ``eigen-character" of $ \sH\left(\wt\RG, K_{\wt \RG}^\circ\right)$ associated to $\pi$.
\end{enumerate}
 
 Moreover, if such a representation $\pi$ is tempered, then 
 $$
 \Int_\pi(f)=2^{-\beta_\pi} \sL'(1/2,\pi) \prod_{v<\infty}\BI_{\pi_v}(f_v).
 $$
\end{conjecture}

 \subsection{The case of regular integral models with parahoric levels}
 We could also state an alternative version of the conjecture using  the arithmetic intersection theory of 
Arakelov and Gillet--Soul\'e \cite[\S4.2.10]{GS} on arithmetic Chow groups. The statement is more straightforward but there is an ambiguity in the choice of Green currents.

In fact, we can introduce a certain parahoric level structure at a finite set $S_{\rm in}$ of {\em inert} places and modify the level structure at a subset $S_{\rm ram}$ of ramified places.  We refer to \cite[\S5]{RSZ3} for the detailed definition for the cases where we have a regular integral model $\CM_{K}\bigl(\wt \RG\bigr)$ (resp.  $\CM_{K_{\wt \RH}}\bigl(\wt \RH\bigr)$) of the Shimura variety $\Sh_{ K}(\wt \RG)$  (resp.  $\Sh_{K_{\wt \RH}}\bigl(\wt \RH\bigr)$) for an appropriate compact open subgroup $K\subset \wt \RG(\BA_f)$ (resp. $K_{\wt \RH}\subset \wt \RH(\BA_f)$).

For a place $\nu: E\incl \BC$, we fix a choice of an admissible Green current $g_{\nu}$ for the (unmodified) cycle $z_K$  on the  complex analytic space $ \CM_{K}\bigl(\wt \RG\bigr)\times_{\Spec E,\nu}\Spec\BC$.  We let $\wt{z}_K$ denote the element in the arithmetic Chow group $\wh\Ch^{n-1}( \CM_{K}\bigl(\wt \RG\bigr))$ defined by the integral model $ \CM_{K_{\wt \RH}} \bigl(\wt \RH\bigr)$ and the Green currents $(g_{\nu})_{\nu\in\Hom_\BQ(E,\BC)}$.   Let $\sH^{\rm spl}(\wt \RG) $ be the spherical Hecke algebra at all {\em split} places of $F'/F$.  Then in \cite[\S4]{RSZ3} we defined a family of Hecke correspondence $\wh R(f)$ on the integral models, for every $f\in \sH^{\rm spl}\left(\wt\RG, K\right)$. Define 
  \begin{align}\label{int f GS}
 \Int_{{\rm GS}} (f)=\Big( \wh R(f)\ast \wh z_{K},\,\,\wh z_{K}\Big)_{{\rm GS}},\quad f\in \sH^{\rm spl}\left(\wt\RG, K\right),
 \end{align}
where ``GS" indicates the arithmetic intersection  theory in Gillet--Soul\'e \cite{GS}.
\begin{conjecture}\label{conj:sph AGGP int}
There exists a choice of admissible Green currents such that there is a decomposition $$
 \Int_{{\rm GS}}(f)=\sum_{\pi}\Int_\pi(f), \quad\text{for all }\, f\in \sH^{ \rm spl}\left(\wt\RG, K\right),
 $$ characterized by the following properties 
  \begin{enumerate}
 \item
  the sum runs over all cohomological (with respect to the trivial coefficient system) automorphic representations of $\wt\RG(\BA)$ with non-trivial $K^S$-invariants for $S=S_{\rm in}\cup S_{\rm ram}$, that are  trivial on $Z^\BQ(\BA)$;
  \item  $\Int_\pi$ is an eigen-distribution for the spherical Hecke algebra $\sH^{\rm spl}(\wt \RG) $ with eigen-character $\lambda_{\pi^{ \rm spl }}$ similarly defined as \eqref{eigen dis}.
  \end{enumerate}
 
 Moreover, if such a representation $\pi$ is tempered and appears in the middle-degree cohomology  $H^{2n-3}(\Sh_{ K_{\wt \RG}}(\wt \RG))$, then 
 $$
 \Int_\pi(f)=2^{-\beta_\pi} \sL'(1/2,\pi) \prod_{v<\infty}\BI_{\pi_v}(f_v).
 $$
 
\end{conjecture}

  \begin{remark} One could formulate a version of the conjecture over a general totally real field, based on \cite[Remark 5.3 and \S8]{RSZ3}. 
  
  \end{remark}
   \begin{remark}\label{rmk:Drin lev}
 In 
 \cite[\S8.2]{RSZ3} we also allowed non-hyperspecial level structure at some {\em split} places and define an integral model using Drinfeld level structure. Though the product Shimura variety may not be regular, we have an ad-hoc definition of the intersection number $\Int(f)$ when $f$ has regular support.
\end{remark}
     \begin{remark}
The recent Thesis of Zhiyu Zhang \cite{ZZ} considers more general maximal parahoric level structure at {\em inert} places than those considered in  \cite{RSZ3}.  It is expected that his result will lead to a version of the above conjecture for such parahoric levels.
 \end{remark}

\subsection{The relative trace formula approach}\label{ss:JR RTF}

Next we move to the relative trace formula approach towards the proof of the conjectures above.  See the articles of
 Beuzart-Plessis \cite{BP22} and Chaudouard \cite{Ch22}  for more background.

We first recall the RTF constructed by Jacquet and Rallis \cite{JR} to attack the Gan--Gross--Prasad conjecture in the Hermitian case.  It is associated to
the triple $(\RG',\RH'_1,\RH'_2)$ where
$$
\RG'=\RR_{F'/F}(\GL_{n-1}\times \GL_{n}),$$
and $$\RH'_1=\RR_{F'/F}\GL_{n-1},\quad \RH'_2=\GL_{n-1}\times\GL_{n},
$$
where $(\RH'_1,\RG')$ is the Rankin--Selberg pair, and $(\RH'_2,\RG')$ the Flicker--Rallis pair.
Moreover it is necessary to insert a quadratic character of $\RH_2'(\BA)$:
$$
\eta=\eta_{n-1,n}:(h_{n-1},h_{n
})\in \RH'_2(\BA)\mapsto \eta^{n-2}_{F'/F}(\det(h_{n-1}))\eta^{n-1}_{F'/F}(\det(h_{n})),
$$
where $\eta_{F'/F}: F^\times\bs \BA^\times\to\{\pm 1\}$ is the quadratic character associated to $F'/F$ by class field theory. 

We introduce (cf.\  \cite[\S3.1]{Z12}) the global distribution on $\RG'(\BA)$ parameterized by a complex variable $s\in\BC$,
\begin{align}
\label{J f'}
\BJ(f',s)=\int_{[\RH'_1]}\int_{[\RH_2']}K_{f'}(h_1,h_2)\,\bigl|\det(h_1)\bigr|^s\,\eta(h_2)\,dh_1\,dh_2,\quad f'\in \mathscr{C}^\infty_c(\RG'(\BA)).
\end{align}
We set its value at $s=0$
\begin{align*}
\BJ(f')=\BJ(f',0),
\end{align*}and its value of the first derivative at $s=0$
 $$
\partial\BJ(f')=\frac d{ds} \Big|_{s=0} \BJ(f',s).
$$
Then the idea is that, in analogy to the usual comparison of two RTFs, we hope to compare the intersection pairing $\Int(f)$ in \eqref{int f} or \eqref{int f GS} and  $\partial \BJ(f')$ for $f\in \sH(\wt\RG, K_{\wt \RG}^\circ)$  and any transfer $f'\in \mathscr{C}^\infty_c(\RG'(\BA))$.

\begin{remark}We comment on the idea behind the construction of Jacquet--Rallis. Note that the Flicker--Rallis period attached to the pair $(\RH'_2,\RG')$ is expected to detect the quadratic base change from the unitary group $\RG$ to $\RG'$ (\cite[ p.144]{F}). Therefore the cuspidal part of the spectral side in ${\rm RTF}_{(\RG',\RH'_1,\RH'_2)}$ is expected to contain only those automorphic representations that are in the image of the quadratic base change from unitary groups. This gives the hope that the spectral sides of the two RTFs should match.
\end{remark}

For $\gamma\in \RG'(F_{v})_\rs$ (here ``$\rs$" indicates ``regular semisimple" in the relative sense, cf. \cite{Z14a}), $f'\in \mathscr{C}^\infty_c(\RG'(F_v))$, and $s\in\BC$, we introduce the (weighted) orbital integral
\begin{align}\label{orb GL}
\Orb(\gamma,f',s)=
\int_{ \RH'_{1,2}(F_v)}f'(h_1^{-1}\gamma h_2)\bigl|\det(h_1)\bigr|^s\eta(h_2)\,dh_1\,dh_2.
\end{align}
We set
\begin{equation}\label{orb del}
   \Orb(\gamma,f') := \Orb(\gamma,f', 0)  \quad\text{and}\quad \del(\gamma,f') := \frac d{ds} \Big|_{s=0} \Orb(\gamma, f',s)  .
\end{equation}

Let $f'=\otimes_v f'_v$ be a pure tensor with {\em regular support}  in the sense that there is a place $u_0$ of $F$ where $f'_{u_0}$ has support in the regular semi-simple locus $\RG'(F_{u_0})_\rs$. 
Then we have a decomposition of \eqref{J f'} into a sum over the set of regular semisimple orbits $ [\RG'(F)_\rs]$,
\begin{align}\label{eq J f' s}
\BJ(f',s)=\sum_{\gamma\in [\RG'(F)_\rs]} \Orb(\gamma,f',s),
\end{align}
where each term is a product of local orbital integrals \eqref{orb GL},
\begin{equation*}
\Orb(\gamma,f',s)=\prod_{v}\Orb(\gamma,f'_v,s).
\end{equation*}

We have a similar definition of orbital integral on the unitary side.  For $g\in \RG(F_v)_\rs$ and $f\in \mathscr{C}^\infty_c(\RG(F_v))$ we introduce the orbital integral
\begin{align}\label{orb U}
\Orb(g,f)=
\int_{ \RH_{1,2}(F_v)}f(h_1^{-1}g h_2)\,dh_1\,dh_2.
\end{align}
We now recall the notion of (smooth) transfer of test functions. At a place $v$, we say that $f\in \mathscr{C}^\infty_c(\RG(F_v))$ and $f'\in \mathscr{C}^\infty_c(\RG'(F_v))$ are transfer of each other if for every $\gamma \in\RG'(F_v)_{\rs}$ we have
\begin{align}\label{orb U}
  \omega(\gamma) \Orb(\gamma,f')=\begin{cases}\Orb(g,f),& \text{if $\gamma$ matches $g\in \RG(F_v)_{\rs}$},\\
   0,&\text{otherwise}.
   \end{cases}
\end{align}Here $\omega(\gamma)$ is a suitable transfer factor. The existence of transfer is known at $p$-adic places by \cite{Z14a} and for many test functions at archimedean places \cite{Xue}.

We can then state an arithmetic intersection conjecture for the arithmetic diagonal cycle on the global integral model $\CM_{K_{\wt \RG}}(\wt \RG)$ \cite[\S8.1]{RSZ3}.
\begin{conjecture}\label{conj I=dJ}
Let $f=\otimes_{v\leq \infty} f_v\in \sH^{\rm spl}(\wt\RG, K_{\wt \RG})$, and let $f'=\otimes_{v\leq \infty} f'_v\in \sC_c^\infty(\RG'(\BA), K')$ be a Gaussian transfer of $f$. Let $S$ denote the finite set of bad places (including all archimedean places). Then there exists a correction function $f'_{\corr,S}\in \sC_c^\infty(\RG'(\prod_{v\in S}F_v))$ depending only on $f_S=\otimes_{v\in S}f_v$ and $f'_S=\otimes_{v\in S}f_v'$, such that 
\begin{align}\label{eq:Int=J}
\Int(f)=-\partial\BJ(f')-\BJ(f'_{\corr,S}\otimes f'^{S})
\end{align}
where $f'^{S}$ denotes $\otimes_{v\notin S} f'_v$.
 Furthermore, for any given $f_S$ we may choose $f'_{S}$ such that $f'_{\corr,S}=0.$
 \end{conjecture}
 Here $f'=\otimes_{v\leq \infty} f'_v\in \sC_c^\infty(\RG'(\BA))$ is a Gaussian transfer of $f=\otimes_{v\nmid\infty}f_{v}\in \sH^{\rm spl}(\wt\RG, K_{\wt \RG})$ if $f_v'$ is a transfer of $f_v$ (resp., the constant function ${\bf 1}$ on the {\em compact} unitary group) for every $v\nmid \infty$ (resp.,  for every $v\mid \infty$).  The existence of Gaussian transfer at $v\mid\infty$ is proved in \cite{BPLZZ}. 
 
 This conjecture implies Conjecture \ref{conj:sph AGGP int} at least for a certain class of Hecke functions $f$. In fact, thanks to the technique of isolating cuspidal spectrum in \cite{BPLZZ}, one can produce many $f$ with transfers $f'$ such that $\BJ(f',s)$ has a spectral decomposition as a sum over cuspidal automorphic representation of $\RG'(\BA)$.  Then the desired assertion follows from the local spherical character identities at all places.

\subsection{Hecke functions with regular supports}
\label{ss:reg}
The comparison  in the equation \eqref{eq:Int=J} can be localized for a large class of test functions $f$ and $f'$. 
Let $f=\otimes_v f_v$ be a pure tensor with {\em regular support} in the sense that there is a place $u_0$ of $F$ where $f_{u_0}$ has support in the regular semisimple locus $\wt\RG(F_{u_0})_\rs$.\footnote{Strictly speaking here we need to bring in non-trivial level structure at split places, see Remark \ref{rmk:Drin lev}. Otherwise we do not know how to produce examples of functions with regular semisimple support.} Then the cycles $\wh R(f)\wh z_{K} $ and $\wh z_{K}$ do not meet in the generic fiber  $\Sh_{K_{\wt \RG}}(\wt \RG)$. The arithmetic intersection pairing then localizes as a sum over all places $w$ of the reflex field $E$ 
$$
\Int(f)=\sum_{w}\Int_w(f).
$$ 
Here for a non-archimedean place $w$, the local intersection pairing $\Int_w(f)$  is defined through the Euler--Poincar\'e characteristic of a derived tensor product on $\CM_{K_{\wt\RG}}(\wt \RG)\otimes_{O_E} O_{E, w}$, cf.\ \cite[4.3.8(iv)]{GS}. 

Similarly, on the GL-side, for $f'=\otimes_v f'_v$ be a pure tensor with {\em regular support}, by \eqref{eq J f' s}
the first derivative $\partial \BJ(f')$ then localizes as  a sum over  places $v$ of $F$,
$$
\partial\BJ(f')=\sum_v\partial\BJ_v(f'),
$$
where   the summand $\partial\BJ_v(f')$ takes the derivative of the local orbital integral (cf.\ \eqref{orb del}) at the place $v$,
\begin{equation}\label{eq: del J v}
\partial\BJ_v(f')= \sum_{\gamma\in [\RG'(F)_\rs]} \del(\gamma,f'_v)\cdot \prod_{u\neq v} \Orb(\gamma,f'_u).
\end{equation}

It is then natural to expect a place-by-place comparison between $\partial\BJ_v(f')$ and $$\Int_v(f):=\sum_{w|v}\Int_w(f)$$ over the places $w$ of $E$ lying above $v$. 

If a place $v$ of $F$ splits in $F'$ (and under the above regularity condition on the support of $f$ and of $f'$), we have \cite[Thm. 1.3]{RSZ3}
$$
\Int_{v}(f)=\partial \BJ_{v}(f')=0.
$$
For every place $w$ of $E$ above a non-split place $v$ of $F$,\footnote{There is also a requirement on  the CM type chosen to define the RSZ moduli functor; let us ignore the subtlety in this notes, and we refer to \cite{RSZ3} for interested readers.} we have a smooth integral model $\CM_{K_{\wt\RG}}(\wt \RG)\otimes_{O_E} O_{E, w}$ when $K_{\wt\RG,v}$ is a hyperspecial compact open subgroup $\wt \RG(O_{F,v})$ (resp.\ a special parahoric subgroup $K^\circ_{\wt\RG,v}$ introduced before) for inert $v$ (resp.\ ramified $v$).  For such places $v$, we have a ``semi-global" conjecture \cite[\S8]{RSZ3}.
\begin{conjecture}\label{conj I=dJ}
Let $f=\otimes_{u\leq \infty} f_u\in \sH^{\rm spl}(\wt\RG, K_{\wt \RG})$, and let $f'=\otimes_{u\leq \infty} f'_u\in \sC_c^\infty(\RG'(\BA), K')$ be a Gaussian transfer of $f$. Then there exists a correction function $f'_{v,\corr}\in \sC_c^\infty(\RG'(F_{v}))$ depending only on $f_v$ and $f'_v$ such that
$$
\Int_{v}(f)=-\partial\BJ_{v}(f')-\BJ(f'_{v,\corr}\otimes f'^{v}),
$$
where $f'^{v}$ denotes $\otimes_{u\neq v} f'_u$. Furthermore, for any given $f_v$ we may choose $f'_v$ such that $f'_{v,\corr}=0.$
 \end{conjecture}

 \begin{theorem}\label{thm:semi-global}
Conjecture \ref{conj I=dJ} holds for an inert $v$ with a hyperspecial level provided the residue characteristic is odd.
\end{theorem}

We sketch the proof.
For $i=n-1,n$, let $W_{i}[v]$ be the pair of nearby Hermitian spaces, i.e., the Hermitian space (with respect to $F'/F$) that is totally negative at archimedean places, and is not equivalent to $W_i$ at $v_0$. Let $\wt\RG[v]$ be the corresponding group, an inner form of $\wt \RG$. We fix an isomorphism $\wt \RG[v_0](\BA^{v}_f)\simeq \wt \RG(\BA^{v}_f)$. Let 
$f=\otimes_u f_u$ be a pure tensor such that 
\begin{altenumerate}
  \item $f_{v}={\bf 1}_{\wt \RG(O_{F,v})}$, and
\item  there is a place $u_0$ of $F$ where $f_{u_0}$ has support in the regular semisimple locus $\wt\RG(F_{u_0})_\rs$.
\end{altenumerate}
Then we have a sum over orbits:
\begin{align}\label{eq A RTF v0}
\Int_{v}(f)=2\log q_v\sum_{g\in \wt\RG[v](F)_\rs}\Int_{v}(g)\cdot\prod_{u\neq v}\Orb(g,f_u),
\end{align}
where the local intersection number $\Int_{v}(g)$ is defined by \eqref{defintprod} in the next section, and $q_v$ is the residue cardinality of $O_F$ at $v$. 
For $F=\BQ$, this is proved \cite[Thm. 3.9]{Z12}. The general case is established in (the proof of) \cite[Thm. 8.15]{RSZ3}.

By the formulas \eqref{eq: del J v} and \eqref{eq A RTF v0},  the comparison between $\Int_{v_0}(f)$ and $\partial\BJ_{v_0}(f')$ is then reduced to a local conjecture that we will consider in the next section, the {\em arithmetic fundamental lemma}:
\begin{align}\label{eq AFL 0}
2\Int_{v}(g)\log q_v=-\omega_v(\gamma) \del(\gamma,{\bf 1}_{\wt \RG(O_{F,v})})
\end{align}
whenever $g$ and $\gamma$ match. (Here $\omega_v(\gamma)\in\{\pm1\}$ is a certain transfer factor.) Now Theorem \ref{thm:semi-global} follows from Theorem \ref{U AFL} in the next section.

\begin{remark}
For a ramified place $v_0$, the analogous question is also reduced to the local {\em arithmetic transfer} conjecture formulated by Rapoport, Smithling and the author in \cite{RSZ1,RSZ2}. We have a result similar to \eqref{eq A RTF v0}, \cite[Thm. 8.15]{RSZ3}.

\end{remark}

\begin{remark}
One may also enlarge the scope of the conjecture by replacing $\sH^{\rm spl}(\wt\RG, K_{\wt \RG})$ by a bigger (spherical) Hecke algebra, e.g., to include the full  spherical Hecke algebra at all the inert places. For the purpose of separating spectrum in $\BJ(f',s)$, Ramakrishnan's density theorem shows that we do not gain more information by doing so. However, besides being a natural and interesting question, there is also another reason to investigate such a question, mainly from the perspective of $p$-adic height pairing (see the next subsection). We hope to pursue this direction in \cite{LRZ2}.

\end{remark}

\subsection{$p$-adic height pairing and a $p$-adic AGGP conjecture}\label{ss:pAGGP}
We now resume the notation from \S\ref{sss std conj}.
Nekov\'a\v{r} defined  a $p$-adic analog of the $\BR$-valued height pairing \eqref{eqn BB}, depending on a suitable splitting of the Hodge filtration:
 \begin{equation}
\label{eqn pHT}
   \xymatrix{\sform_{\it p}\colon \Ch^i(X)_{\BQ,0}\times \Ch^{j}(X)_{\BQ,0}\ar[r]& \BQ_p,\quad i+j-1=d=\dim X},
\end{equation}
which is conditional on more accessible conjectures than  needed to define \eqref{eqn BB} (essentially the weight-monodromy conjecture at all places), see \cite{Ne1}.  In fact assuming the  weight-monodromy conjecture (including the $p$-adic analog) from now on,  the $p$-adic Abel--Jacobi map lands in the Bloch--Kato Selmer group (a finite dimensional $\BQ_p$-vector space)   \cite{Ne2}:
$$
\xymatrix{\Ch^i(X)_0\ar[r]& H^1_f(E,H^{2i-1}(X_{\ov E},\BQ_p(i))).}
$$
The $p$-adic height pairing \eqref{eqn pHT} then factors through a pairing on Selmer groups
 \begin{equation}
\label{eqn pHT Sel}
 \xymatrix{ \sform_{{\it p}}\colon H^1_f(E,H^{2i-1}(X_{\ov E},\BQ_p(i)))\times  H^1_f(E,H^{2j-1}(X_{\ov E},\BQ_p(i)))\ar[r]& \BQ_p}.
\end{equation}
Therefore the $p$-adic height pairing depends only on the absolute cohomological class of cycles, and hence it is more accessible than the $\BR$-valued height pairing \eqref{eqn BB}.  

Back to the setting of AGGP, we assume that  the middle degree cohomology  decomposes as a $ \sH\left(\RG, K\right)\times \Gal(\ov F/F')$-module (see \cite{KSZ}
for recent progress towards this problem)
\begin{align}\label{eq:decomp H}
H^{2n-3}(\Sh_K(G)_{\ov F},\ov\BQ_p(n-1))=\bigoplus_{\pi} \pi_f^K\boxtimes \rho_{\pi}^\vee
\end{align}
where $\pi$ runs over all cohomological automorphic representations of $G$ and $\rho_{\pi}$ is a representation of $\Gal(\ov F/F')$ with $\ov\BQ_p$-coefficient, which up to semisimplification is the (normalized) global Langlands correspondence of $\pi$, at least for stably cuspidal $\pi$ (note that $G$ is a product of two unitary groups and $\rho_{\pi}$ is understood as the tensor product of two Galois representations).  Here we fix an isomorphism  $\ov\BQ_p\simeq\BC$ to transport $\BC$-valued automorphic representations to $\ov\BQ_p$-valued ones. For $\varphi^\vee\in (\pi_f^\vee)^K\simeq (\pi_f^K)^\ast$ (here $\pi_f^\vee$ denotes the contragredient and $(\pi_f^K)^\ast$ denotes the usual linear dual),  we have a linear map  $\varphi^\vee(\cdot): \pi_f^K\to \ov\BQ_p $ and an induced linear map
$$
\xymatrix{H^1_f(F',H^{2n-3}(\Sh_K(G)_{\ov F},\ov\BQ_p(n-1)))\ar[r]^-\sim& \bigoplus_{\pi} \pi_f^K\otimes H^1_f(F',\rho_{\pi}^\vee)\ar[d]^-{\varphi^\vee(\cdot)\otimes \id} \\
 &H^1_f(F',\rho_{\pi}^\vee) .}
$$
In particular, we apply the construction to the class of the modified arithmetic diagonal cycle $z_{K,0}$ and obtain an element  which we denote by
\begin{align}\label{eq: sel cl}
\int_{z_{K,0}}\varphi^\vee \in H^1_f(F',\rho_{\pi}^\vee).
\end{align}
As the notation suggests, this may be viewed as an analog of the automorphic period integral \eqref{eq:P pi}. In  a work in progress of Daniel Disegni and the author \cite{DZ}, we formulate  a $p$-adic AGGP conjecture. The refined form of this conjecture is an exact identity relating the $p$-adic height pairing 
$$
\left(\int_{z_{K,0}}\varphi ,\int_{z_{K,0}}\varphi^\vee\right)_{\it p} , \quad \varphi\in \pi^K,  \, \varphi^\vee\in (\pi^\vee)^K
$$ for tempered cuspidal $\pi$, ordinary at all $p$-adic places,
to the first derivative of a suitable $p$-adic L-function  $\sL'_{\it p}(s,\pi)$ times the local canonical invariant forms given by Ichino--Ikeda \eqref{alpha1},
$$
2^{-\beta_\pi} e_p(\pi)^{-1} \sL'_{\it p}(1/2,\pi) \prod_{v<\infty}\alpha_{v}(\varphi_v,\varphi^\vee_v),
$$possibly up to certain local modification $e_p(\pi)^{-1}$ at places above $p$.  In particular, this formula would imply that the class $\int_{z_{K,0}}\varphi^\vee \in H^1_f(F,\rho_{\pi}^\vee)$ in the Bloch--Kato Selmer group does not vanish for some $\varphi\in \pi_f$ if and only if  $\ord\sL_{\it p}(s,\pi)=1$ and the  local canonical invariant forms do not vanish.  

In \cite{DZ} we prove this conjectural formula when all the $p$-adic places are split in $F'/F$, at least for $\pi$ with mild ramifications and with good ordinary reduction at all $p$-adic places. The overall strategy is similar to the relative trace formula approach: we study  the $p$-adic analog of \eqref{eqn BB K}
  \begin{align}\label{p int f}
 \Int_{{\it p}}(f)=\Big( R(f)\ast z_{K,0},\,\,z_{K,0}\Big)_{{\it p}},\quad f\in \sH\left(\RG, K\right)
 \end{align}
 and a $p$-adic interpolation of the Jacquet--Rallis RTF in \S\ref{ss:JR RTF}. Assume that the cycles have disjoint support. Then the $p$-adic height pairing also localizes into a sum of local heights over all the (non-archimedean) places. For almost all places away from $p$, the local pairing is essentially the same as the intersection pairing in Theorem \ref{thm:semi-global} (except replacing $\log q_v$ in \eqref{eq A RTF v0} by a $p$-adic $\log_p q_v$). For the remaining places away from (hyperspecial) $p$ and mostly the places above $p$, substantial work is needed and is beyond the scope of this notes. In some way, the difficulty for the local $p$-adic height at $p$-adic places resembles that for the local $\BR$-valued height \eqref{eqn BB} at archimedean places. Under the ordinary (at $p$-adic places) hypothesis, we are able to show a suitable limit of the local heights at $p$-adic places vanishes, resembling the treatment by Perrin-Riou of the $p$-adic Gross--Zagier formula.

\section{AFL and AT conjecture for local Shimura varieties}
\label{s:loc}
In this section we present more detail on the two local questions, namely the arithmetic analogs of the fundamental lemma and of (smooth) transfer, which are called the arithmetic fundamental lemma (AFL) and arithmetic transfer (AT). 
\subsection{Unitary Rapoport--Zink spaces}\label{ss:RZ}

We first recall the unitary Rapoport--Zink (RZ) spaces, a local analog of CM cycles, local KR divisors, and then we  state the AFL conjecture, cf. \cite[\S4]{RSZ2} and \cite[\S3]{Z21}.

Let $F'/F$ be an unramified quadratic extension of $p$-adic local fields with $p$ {\em odd} and let $n\geq 1$. Denote by $q$ the residue cardinality of $O_F$. 
Let $\breve F$ denote the completion of a maximal unramified extension of $F$. For $\Spf O_{\breve F}$-schemes $S$ (i.e. a $O_{\breve F}$-scheme on which $p$ is locally nilpotent), we consider triples $(X, \iota, \lambda)$, where
\begin{enumerate}
\item
$X$ is a $p$-divisible group of absolute height $2nd$ and dimension $n$ over $S$, where  $d:=[F: \BQ_p]$, 
\item  $\iota$ is an action of $O_{F'}$ such that the induced action of $O_{F}$ on $\Lie X$ is via the structure morphism $O_{F}\to \CO_S$,  and
\item $\lambda: X\to X^\vee$ is a principal ($O_{F}$-relative) polarization. 
\end{enumerate}
Hence $(X, \iota|_{O_{F}})$ is a strict $O_{F}$-module of relative height $2n$ and dimension $n$. We require that the Rosati involution $\Ros_\lambda$ induces on $O_{F}$  the non-trivial Galois automorphism in $\Gal(F'/F)$, denoted by $ O_{F'}\ni a\mapsto \ov a$, and that the \emph{Kottwitz condition} of signature $(n-1,1)$ is satisfied, i.e.
\begin{equation}\label{kottwitzcond}
   \charac \bigl(\iota(a)\mid \Lie X;\, T\bigr) = (T-a)^{n-1}(T-\ov a) \in \CO_S[T]
	\quad\text{for all}\quad
	a\in O_{F'} . 
\end{equation} 
An isomorphism $(X, \iota, \lambda) \isoarrow (X', \iota', \lambda')$ between two such triples is an $O_{F'}$-linear isomorphism $\varphi\colon X\isoarrow X'$ such that $\varphi^*(\lambda')=\lambda$.

Over the residue field $\mbF$ of $O_{\breve {F}}$, there is a triple $(\BX_n, \iota_{\BX_ n}, \lambda_{\BX_n})$ such that $\BX_n$ is  supersingular, unique up to $O_{F'}$-linear quasi-isogeny compatible with the polarization. We fix such a triple which we call the {\em framing object}. Then there is the Rapoport--Zink formal moduli scheme $\CN_n = \CN_{n, F'/F}$ associated to the unitary group for the quasi-split 
$n$-dimensional Hermitian $F'$-vector space.
 By definition  $\CN_n$ represents the functor over $\Spf O_{\breve F}$ that associates to each $S$ the set of isomorphism classes of quadruples $(X, \iota, \lambda, \rho)$ over $S$, where the last entry is an $O_F$-linear quasi-isogeny of height zero defined over the special fiber $\ov S:=S\times_{ \Spf O_{\breve F}} \Spec \mbF$,
\[
   \rho \colon X\times_S\ov S \to \BX_n \times_{\Spec \mbF} \ov S,
\]
such that $\rho^*((\lambda_{\BX_n})_{\ov S}) = \lambda_{\ov S}$. The map $\rho$ is called the {\em framing}.
The formal scheme $\CN_n$ is  formally locally of finite type and formally smooth over $\Spf O_{\breve {F}}$ of relative dimension $n-1$. 

The group of quasi-automorphisms of the framing object is
$$\Aut^\circ (\BX_n,\iota_{\BX_n},\lambda_{\BX_n}) = \{g\in \End^\circ_F(\mbX_n)\mid\ g^\vee\circ \lambda_{\mbX_n} \circ g = \lambda_{\mbX_n}\}.$$
The condition $g^\vee \circ \lambda_{\mbX_n}\circ g = \lambda_{\mbX_n}$ may also be formulated as $gg^* = \mr{id}$, where $g \mapsto g^\ast=\Ros_{\lambda_{\BX_n}}(g)$ denotes the Rosati involution induced by the polarization $ \lambda_{\mbX_n}$. Then $\Aut^\circ (\BX_n,\iota_{\BX_n},\lambda_{\BX_n})$ acts on $\CN_n$ by changing the framing: 
$$g\cdot (X,\iota,\lambda,\rho) = (X,\iota,\lambda, g \circ \rho).$$

Another description is as follows. Taking $n=1$, there is a unique triple $(\mbE, \iota_{\mbE}, \lambda_{\mbE})$ over $\mbF$ with signature $(1,0)$. Set
\begin{align}\label{eq:BVn}
\BV_n:=\Hom^\circ_{O_{F'}}(\BE, \BX_n),
\end{align}called the space of {\em special homomorphisms}. It has the structure of an $n$-dimensional Hermitian $F'$-vector space with respect to the Hermitian form
$$
\pair{x,y}= \lambda_{\mbE}^{-1}\circ y^\vee\circ \lambda_{\mbX_n} \circ x\in \End^\circ _{F'}(\mbE)\simeq F'.
$$
It is the unique (up to isomorphism) $n$-dimensional Hermitian space that does {\em not} contain a self-dual $O_F$-lattice.
 Then there is a natural isomorphism
\begin{equation}\label{Aut cong U}
   \Aut^\circ (\BX_n,\iota_{\BX_n},\lambda_{\BX_n}) \cong \U\bigl(\BV_n\bigr)(F), 
\end{equation}
where $\Aut^\circ (\BX_n,\iota_{\BX_n},\lambda_{\BX_n})$ acts by composition on $\mbV_n$.

We also define a broader class of RZ spaces $\CN_n^{[t]}$ for an integer $0\leq t\leq n$. We only need to replace the principal polarization in the triple $(X,\iota,\lambda)$ by the following variant:
\begin{enumerate}
\item[(3)] $\lambda: X\to X^\vee$ is a polarization such that $\ker(\lambda)\subset X[\varpi]$ has height $q^{2t}$. 
\end{enumerate}
In particular, we have $\CN_n=\CN_n^{[0]}$. 
The resulting RZ space $\CN_n^{[t]}$ is formally locally of finite type and of strictly semi-stable reduction over $\Spf O_{\breve {F}}$ with relative dimension $n-1$ \cite{Cho,Go}. The case of $t=1$ is of special interest, termed as the ``almost unramified level" in \cite{RSZ2}, see also \S\ref{ss:ATC} below. In general $\CN_n^{[t]}$ is attached to a parahoric subgroup fixing a vertex lattice of type $t$ in a Hermitian space $V$ whose discriminant has valuation $\equiv t\mod 2$.  Then the space $\BV_n$ (depending on $t$ but we suppress $t$ from the notation) of special homomorphisms defined by \eqref{eq:BVn} has the structure of a Hermitian space of  discriminant with valuation $\equiv 1+t\mod 2$. Moreover, there is a natural isomorphism 
$\CN_n^{[t]}\simeq \CN_n^{[n-t]}$ \cite[\S5.1]{ZZ}. In particular we have $\CN_1^{[0]}\simeq \CN_1^{[1]}$. When $n=2$, $\CN_2^{[0]}\simeq \CN_2^{[2]}$ is isomorphic to the  Lubin--Tate space of height two, and $\CN_2^{[1]}$ is isomorphic to Deligne's semistable integral model of the Drinfeld half plane $\wh\Omega_{2}$ \cite{KR14}.

\subsection{The AFL conjecture/theorem}\label{ss:AFL}

Let $\CN_{n-1,n}=\CN_{n-1}\times_{\Spf O_{\breve F}}\CN_n$. Then $\CN_{n-1,n}$ admits an action by $G(F)$ where we denote now $G=\RU(\BV_{n-1})\times \RU(\BV_{n})$.

There is a natural closed embedding $\delta \colon \CN_{n-1}\to \CN_n$  (a local analog of the closed embedding \eqref{inc M}), see \eqref{em} below for more detail. Let 
\begin{equation*}
   \xymatrix{\Delta \colon \CN_{n-1} \xra{\,}   \CN_{n-1,n}}
\end{equation*}
be the graph morphism of $\delta$. We denote by $\Delta_{\CN_{n-1}}$ the image of $\Delta$. For $g\in \RG(F)_\rs$, we consider the intersection number on $\CN_{n-1, n}$ of $\Delta_{\CN_{n-1}}$ with its translate $g\Delta_{\CN_{n-1}}$, defined through the derived tensor product of the structure sheaves (see also the beginning of \S\ref{ss:reg}),
\begin{equation}\label{defintprod}
   \Int(g) := \left( \Delta_{\CN_{n-1}}, g\cdot\Delta_{\CN_{n-1}}\right)_{\CN_{n-1, n}} := \chi\left({\CN_{n-1, n}},  \CO_{\Delta_{\CN_{n-1}}}\Ltimes_{\CO_{\CN_{n-1, n}}}\CO_{g\cdot\Delta_{\CN_{n-1}}}\right) . 
\end{equation}Here, for a finite complex $\CF$ of coherent $\mcO_{\CN_{n-1,n}}$-modules, we define its Euler--Poincar\'e characteristic as 
$$
\chi(\CN_{n-1,n},\CF)=\sum_{i,j} (-1)^{i+j} \mr{len}_{O_{\breve F}} H^{j}(\CN_{n-1,n},H_i(\CF))
$$
if the lengths are all finite. The regularity of $g$ implies that the (formal) schematic intersection $\Delta_{\CN_{n-1}}\cap g\cdot \Delta_{\CN_{n-1}}$ is a {\em proper scheme} over $\Spf O_{\breve{F}}$, and hence the finiteness of the Euler–Poincar\'e characteristic.

Recall from  \eqref{orb del} the derivative of the orbital integral. 
\begin{theorem}[Arithmetic Fundamental Lemma (AFL)]
\label{U AFL}
Let $\gamma\in\RG'(F)_\rs$ match an element $g\in \RG(F)_\rs$. Then  
\begin{equation*}\label{introAT}
\omega(\gamma)   \del\left(\gamma,{\bf 1}_{\RG'(O_{F})}\right)=- 2\,\Int(g)\cdot\log q.
\end{equation*}
Here $\omega(\gamma) \in\{\pm 1\}$ is an explicit transfer factor (see \cite{RSZ2}).
\end{theorem}
For $F=\BQ_p$, the theorem is proved \cite{Z21} in large residue characteristic. The general case is established in  \cite{MZ} in large (meaning $p\geq n$) residue characteristic and in \cite{ZZ} in any  {\em odd} residue characteristic.  A survey on the conjecture can be found in \cite{Z12b}, and a survey on the proof of AFL can be found in \cite{Z19}; see also Remark \ref{rem AFL}.

This statement was conjectured in \cite{Z12}, where the low rank case of $n=2,3$ was proved (a simplified proof was given by Mihatsch in \cite{M-AFL}).
 Rapoport, Terstiege, and the author in \cite{RTZ} proved a special case, the so-called {\em minuscule} case, for $p$ large comparing to $n$  (a simplified proof was given by Li and Zhu in \cite{LZ}).

\begin{remark}
Note that in the formulation of AFL we have restricted ourselves to odd residue characteristic $p$. This is caused by the same assumption in the theory of RZ spaces \cite{RZ}. We expect the statement in the AFL conjecture to also hold when $p=2$.
\end{remark}

We refer to \cite[\S4]{RSZ2} for some other equivalent formulations of the AFL conjectures. Moreover, in \cite{Z23} we also formulated an AFL conjecture for the general Bessel cycles, generalizing the arithmetic diagonal cycle $\Delta_{\CN_{n-1}}$. 
However, there are no global analogs of these local cycles due to the lack of relevant Shimura varieties.  The corresponding FL conjecture was formulated by Yifeng Liu in \cite{Liu14}, which is still unproven.

\subsection{The AFL conjecture of Liu in the Fourier--Jacobi $\RU(n)\times\RU(n)$ case}
\label{ss:Lie AFL}

A closely related statement is the AFL conjecture in the Fourier--Jacobi $\RU(n)\times\RU(n)$ case, formulated by Yifeng Liu \cite{Liu21}. This statement is called the semi-Lie version of AFL conjecture in \cite{Z21,MZ,ZZ} and is shown to be equivalent to the AFL (the so-called ``group version") in Theorem \ref{U AFL}. In fact the interplay between these two statements is crucial to the inductive proof  of both simutanously!

 In \cite{KR-U1}, Kudla and Rapoport have defined for every non-zero $u\in \BV_n$ a special divisor $\CZ(u)$ on $\CN_n$, which for simplicity will be called a KR divisor. For its definition, note that $\mcN_1 \iso \Spf O_{\breve F}$, so $(\mbE, \iota_{\mbE}, \lambda_{\mbE})$ deforms to a unique triple $(\mcE, \iota_{\mcE}, \lambda_{\mcE})$ over $O_{\breve F}$, called its \emph{canonical lift}. In other words, this (with the $\rho_{\mcE}$ for the deformation) is the universal object over $\CN_1$ with Galois conjugated $O_{F'}$-action. Then the KR-divisor $\mcZ(u)$ is defined as the locus where the quasi-homomorphism $u\colon \BE\to \BX_n$ lifts to a homomorphism from $\CE$ to the universal object over $\CN_n$. By \cite[Prop.\ 3.5]{KR-U1}, $\CZ(u)$ is a relative Cartier divisor (or empty). It follows from the definition that if $g\in \U(\BV_n)(F)$, then
\begin{align}\label{act KR}
g \CZ(u)=\CZ(gu).
\end{align}

The KR divisor plays a fundamental role in the formulation of the arithmetic Siegel--Weil formula, where one is interested in the intersection of $n$ KR divisors, see \eqref{eq KR conj} below. It also appears naturally in Liu's formulation of the AFL conjecture in the Fourier--Jacobi case, as we explain now. 
For simplicity we will write $\CN_n\times\CN_n$ for the fiber product $\CN_n\times_{\Spf O_{\breve F}}\CN_n$.
For $g\in \U(\BV_n)(F)$, let $\Gamma_g\subseteq \CN_n\times\CN_n$ be the graph of the automorphism of $ \CN_n$ induced by $g$.
The fixed point locus of $g$ is defined as the intersection
\begin{align}\label{Ng}
\CN_n^g :=\Gamma_g\,\cap \Delta_{\CN_n},
\end{align}
viewed as a closed formal subscheme of $\CN_n$.
We also form the {\em  derived fixed point locus}, denoted by $\LN^g_n$, i.e. the derived tensor product
\begin{align}\label{der Ng}
\LN^g_n :=\Gamma_g\,\jiao \Delta_{\CN_n}:=\CO_{\Gamma_g}\Ltimes_{\CO_{\CN_n\times\CN_n}}\CO_{\Delta_{\CN_n}}
\end{align}
viewed as an element in the K-group $K_0^{\CN_n^g}(\CN_n)$ with support, cf. \cite[Appendix B]{Z21}. 

For a pair $(g,u)\in  (\U(\BV_n)\times  \BV_n)(F)_\rs$ (here ``$\rs$" indicates ``regular semisimple" in the relative sense, i.e., $\{g^i u\}_{i=0}^{n-1}$ form a basis of  the vector space $\BV_n$), we now set
\begin{equation}\label{def int g u}
   \Int(g,u) := \chi\left({\CN_{ n}}, \, \mcO_{\mcZ(u)} \Ltimes_{\CO_{\CN_n}}\, \LN_n^g\right) . 
\end{equation}
When $(g,u)$ is regular semi-simple, the intersection $\CZ(u)\cap\CN^g_n$ is a proper {\em scheme} over $\Spf O_{\breve{F}}$ and hence the right-hand side of \eqref{def int g u} is finite.
The number $\Int(g,u)$ depends only on the $\U(\BV_n)(F)$-orbit of $(g,u)\in ( \U(\BV_{n})\times  \BV_{n})(F)$. There is an equivalent definition that does not involve the derived fixed point locus  $\LN^g_n$ (cf. \cite[Rem.~ 3.1]{Z21}),
\begin{equation}\label{eq:Int alt}
\Int(g,u)= \chi\left({\CN_{ n}\times\CN_n}, \, \CO_{\Gamma_g}\Ltimes_{\CO_{\CN_n\times\CN_n}} \CO_{\Delta_{\CZ(u)}}\right) . 
\end{equation}

To define the analytic side, we consider the symmetric space
\begin{equation}\label{Sn def}
     S_{n} := \{\,\gamma\in \RR_{F'/F}\GL_n\mid \gamma\ov {\gamma}=1_n\,\},
\end{equation}
and the $F$-vector space
\begin{equation}\label{V' def}
V'_{n}=F^{n}\times (F^{n})^\ast,
\end{equation} 
where $(F^{n})^\ast=\Hom_{F}(F^n,F)$ denotes the $F$-linear dual space. For convenience we will identify $F^{n}$ (resp. $(F^{n})^\ast$) with the space of column vectors (resp. row vectors). With the tautological pairing,
we will view $V_n'$ as an $(F\times F)/F$-hermitian space.
Let $\GL_{n,F}$
act  (diagonally)  on the product $S_{n}\times V'_{n}$,
$$
h\cdot(\gamma, (u_1,u_2))=(h^{-1} \gamma h, ( h^{-1}u_1,u_2 h)).
$$

For $(\gamma, u') \in (S_{n} \times V_{n}')(F)_\rs$, $\Phi' \in C_c^\infty((S_{n} \times
V_{n}' )(F))$, and $s \in \BC$, we define the orbital integral
\begin{equation}\label{Orb(gamma,f',s)}
   \Orb((\gamma,u'),\Phi',s) := \int_{\GL_{n}(F)} \Phi'(h\cdot (\gamma, u')) \lv \det h \rv^s \eta(h)\, dh,
\end{equation}
with special values
\begin{align}\label{eq:orb s=0}
   \Orb((\gamma,u'), \Phi') &:=\Orb((\gamma,u'), \Phi', 0),\ \ \text{and}\\
	\del((\gamma,u'), \Phi') &:=\frac{d}{ds} \Big|_{s=0} \Orb((\gamma,u'), \Phi', s).
\end{align}

\begin{theorem}[AFL, semi-Lie algebra version]\label{AFLconj rs}
\label{AFL lie}
Suppose that $(\gamma,u')\in (S_{n}\times V_n')({F})_\rs$ matches the element $(g,u)\in ( \U(\BV_{n})\times  \BV_{n})(F)_\rs$. Then 
\[
\omega(\gamma,u')\del\bigl((\gamma,u'), \mathbf{1}_{(S_{n}\times V_n')(O_{F})}\bigr) 
	   = -\Int(g,u)\cdot\log q. 
\]Here, $\omega(\gamma,u')\in \{\pm 1\}$ is the transfer factor \cite[(2.13)]{Z21}. 
\end{theorem}
Though our formulation  appears to differ from that of Liu in \cite{Liu21}, they are essentially equivalent.
Note that the factor 2 in Theorem \ref{U AFL} disappears here due to the change of the orbital integrals in the LHS.

The case of a fixed $u$ with $(u,u)\in O_{F}^\times$ is equivalent to the AFL conjecture in Theorem \ref{U AFL}, cf. also \cite[Conj. 3.2]{Z21}. In that case there is a natural isomorphism 
$$
\CZ(u)\simeq \CN_{n-1}.
$$
By \cite[Prop. 4.12]{Z21},
the intersection number \eqref{eq:Int alt} is equal to the one considered in \cite{Z12},
$$
\chi\left({\CN_{ n-1,n}}, \, \CO_{(1,g)\cdot\Delta_{\CN_{n-1}}  }\Ltimes_{\CO_{\CN_{n-1,n}}} \CO_{\Delta_{\CN_{n-1}}}\right),
$$and
the orbital integral reduces to the one in {\it loc. cit.} as well, and hence Theorem \ref{AFLconj rs} implies  Theorem \ref{U AFL}. It is less trivial to show in \cite{Z21,MZ,ZZ} that  Theorem \ref{U AFL} with $n+1$ replacing $n$ implies  Theorem \ref{AFLconj rs}.

\subsection{Arithmetic Transfer Conjectures}\label{ss:ATC}
We finally discuss analogs of the smooth transfer in the arithmetic context over a $p$-adic field. Let $F'/F$ be as above. Fix an integer $0\leq t\leq n-1$ and $e\in\{0,1\}$. Then there is a natural embedding $\delta \colon \CN^{[t]}_{n-1}\to \CN_n^{[t+e]}$. In terms of the moduli functors, $\delta$ sends a quadruple 
\begin{align}\label{em}
( X, \iota_{X}, \lambda_{X},\rho_X)\in \CN^{[t]}_{n-1}(S)\mapsto
  ( X\times \CE, \iota_{X}\times  \iota_{\CE}, \lambda_{X}\times  \lambda_\CE,\rho_X\times \varpi^e\rho_\CE)
  \end{align}
   where $( \CE,\iota_{\CE},  \lambda_\CE, \rho_\CE) $ is the canonical lift.
Let $
\Delta \colon \CN^{[t]}_{n-1} \to  \CN_{n-1,n}$ be the graph morphism of $\delta$. 

Consider the product $\CN_{n-1,n}=\CN^{[t]}_{n-1}\times_{\Spf O_{\breve F}}\CN_n^{[t+e]}$. It is regular precisely when (at least) one of the two 
factors is formally smooth over $\Spf O_{\breve{F}}$; these cases were considered in \cite{RSZ2}. In general, Zhiyu Zhang resolves the singularities  of $\CN_{n-1,n}$ by blowing up a certain Weil divisor of $\CN_{n-1,n}$ in \cite{ZZ}. The resulting space is regular and denoted by $\wt{\CN}_{n-1,n}$. The action of $\RG(F)=(\RU(\BV_{n-1})\times \RU(\BV_{n}))(F)$ lifts to one on $\wt{\CN}_{n-1,n}$.  Moreover, taking the strict transform does not change $\Delta_{ \CN^{[t]}_{n-1}}$ and hence the morphism $\Delta$ lifts to $\wt\Delta: \CN^{[t]}_{n-1} \to \wt \CN_{n-1,n}$. Similar to \eqref{defintprod} we may then define for $g\in \RG(F)_\rs$,

\begin{equation}
   \Int(g) :=\chi\left({\wt\CN_{n-1, n}},  \CO_{\wt\Delta_{\CN^{[t]}_{n-1}}}\Ltimes_{\CO_{\wt\CN_{n-1, n}}}\CO_{g\cdot\wt\Delta_{\CN^{[t]}_{n-1}}}\right) . 
\end{equation}

Let $\Lambda^\flat$ be a Hermitian lattice of type $t$ (i.e., $\Lambda^\flat\subset (\Lambda^\flat)^\vee$ is of colength $t$, where  $(\Lambda^\flat)^\vee$ is the dual lattice of $\Lambda^\flat$). It corresponds to a vertex on the Bruhat--Tits building of the unitary group $\RU(V^\flat)(F)$, where $V^\flat=\Lambda^\flat\otimes_{O_F} F$, and the stabilizer of this vertex is a maximal parahoric  subgroup. Let $\Lambda=\Lambda^\flat\oplus O_{F'} u$ be an orthogonal sum for a vector $u$ whose norm $(u,u)$ has valuation $e\in\{0,1\}$. Then $\Lambda$ is a  Hermitian lattice of type $t+e$. Choose an orthogonal basis $\{u_1,\cdots, u_{n-1}\}$ of $\Lambda^\flat$ such that ${\rm val} (u_i,u_i)\in\{0,1\}$. Using this basis to identify $\Lambda^\flat\simeq O_{F'}^{n-1}\subset F'^{n-1}$. The lattice chain 
$\{\Lambda^\flat,(\Lambda^{\flat})^\vee \}$ defines a parahoric integral model of ${\rm R}_{F'/F}\GL_{n-1}$ and let ${\rm R}_{F'/F}\GL_{n-1}(O_F)$ be the parahoric subgroup. Similarly extend the basis to a basis of $\Lambda$ by appending $u$ and we have a  parahoric subgroup ${\rm R}_{F'/F}\GL_{n}(O_F)$ of ${\rm R}_{F'/F}\GL_{n-1}(F)$. Note that we have a ``good-position" condition ${\rm R}_{F'/F}\GL_{n-1}(O_F)\subset {\rm R}_{F'/F}\GL_{n}(O_F)$.  Let $\RG'(O_F)$ denotes the product of the two parahoric subgroups. 

In  \cite{ZZ} Zhiyu Zhang proved: 
\begin{theorem}[Arithmetic Transfer Conjecture (ATC) at maximal parahoric levels]
\label{U ATC}Assume that $F$ is unramified over $\BQ_p$ when $0<t<n-1$.
Let $\gamma\in\RG'(F)_\rs$ match an element $g\in \RG(F)_\rs$. Then  
\begin{equation*}\label{introAT}
\omega(\gamma)   \del\left(\gamma,c\cdot {\bf 1}_{\RG'(O_{F})}\right)=- 2\,\Int(g)\cdot\log q.
\end{equation*}
Here $\omega(\gamma) \in\{\pm 1\}$ is an explicit transfer factor (see \cite{ZZ}), and $c$ is an explicit constant (depending on $n,t,e$).
\end{theorem}

 \begin{remark}\begin{altenumerate}
  \item In a work in progress \cite{LRZ2} Li, Rapoport, and the authors are formulating more AT conjectures, including new examples involving the space $\CN^{[1]}_{n-1}\times_{\Spf O_{\breve F}}\CN_n^{[0]}$ which does not appear in the above theorem of Zhiyu Zhang.
 \item When $F'/F$ is a {\em ramified} quadratic extension of $p$-adic fields (for odd $p$), Rapoport, Smithling, and the author proposed an  AT  conjecture in  \cite{RSZ1,RSZ2}. We proved the conjecture for $n=2,3$. However, the higher rank case remains open.
 \item A major difficulty to formulate more general arithmetic transfer conjectures is the lack of  the regular integral models  for deeper (than parahoric) level structure.
 \end{altenumerate}
 \end{remark}

\begin{remark}\label{rem AFL}
The strategy of the proof in \cite{Z21,MZ,ZZ} is of global nature, crucially relying on the modularity of generating series of special divisors. In fact, in \cite{Z21} the author uses a theorem of Bruinier--Howard--Kudla--Rapoport--Yang \cite{BHKRY} over the integral models of Shimura varieties, and in \cite{MZ,ZZ} the authors only need a theorem of Yuan--Zhang--Zhang in \cite{YZZ1} over the generic fibers of Shimura varieties. On the other hand, Beuzart-Plessis has found a local proof of the Jacquet--Rallis fundamental lemma \cite{BP21}. Is there a local proof of the AFL conjecture? 
\end{remark}

\section{Some other related works}
\label{s:other}

\subsection{Some other generalizations of Gross--Zagier formula}
Besides the orthogonal and unitary AGGP conjecture in \S\ref{ss: AGGP},
there are several other conjectural generalizations of the Gross--Zagier formula, corresponding to different constructions of algebraic cycles (and possibly different families of L-functions).  Shouwu Zhang summarized some of them in his notes \cite{ZSW10} and since then there have been a couple of new additions.

The first one is the arithmetic Rallis inner product for both orthogonal and unitary groups, which involves the standard L-function of an automorphic representation of $\SO(2n+1)$ and $\RU(n)$. By the doubling method, this is largely tied with 
Kudla's program on the arithmetic analog of the Siegel--Weil formula.  We refer to Chao Li's articles \cite{Li-BAMS,Li-IHES} for more detail on the history and the recent progress. A milestone is the precise formulation of the conjecture due to  Kudla and Rapoport in \cite{KR-U1}, which relates  the intersection numbers of $n$ KR divisors defined in \S\ref{ss:Lie AFL} on $\CN_n$,
\begin{equation}\label{eq KR conj}
   \Int(x_1,\cdots, x_n) := \chi\left({\CN_{ n}}, \, \mcO_{\mcZ(x_1)} \Ltimes \cdots \Ltimes  \mcO_{\mcZ(x_n)} \right),\quad x_1,\cdots, x_n\in \BV_n  
   \end{equation}
to the central derivative of a local density function. The role of their conjecture in the arithmetic Siegel--Weil formula is similar to that of the AFL in the AGGP Setting. Chao Li and the author have proved the Kudla--Rapoport  conjecture in both the unitary and the orthogonal case \cite{LZ1,LZ2}. Together with other non-trivial inputs (particularly a proof of the extension of the Kudla--Rapoport conjecture to the ramified quadratic case), Chao Li
and Yifeng Liu have established the first generalization of  the Gross--Zagier formula  to arbitrary high dimensional Shimura varieties \cite{LL1,LL2}. Similar to the $p$-adic AGGP conjecture in \S\ref{ss:pAGGP}, there is also a $p$-adic arithmetic inner product formula, proved by Disegni and Yifeng Liu under certain constraints on the local ramifications \cite{DL}.

Another generalization is formulated by Yifeng Liu \cite{Liu21}, which is the arithmetic analog of the $\RU(n)\times \RU(n)$  Fourier--Jacobi case of GGP conjecture. As already mentioned in \S\ref{ss:Lie AFL}, the local questions are closely related to the  $\RU(n)\times \RU(n+1)$-case. Note that there is no known conjectural generalization of  the Gross--Zagier formula for the Rankin--Selberg L-function for other pairs $\RU(m)\times \RU(n)$ (a very curious coincidence: the Rankin--Selberg theory takes the simplest form for $\GL(n)\times \GL(m)$ when $|n-m|\in \{0,1\}$). Moreover, we do not know whether Liu's construction has an analog for orthogonal groups. In fact, all the known ``arithmetic conjectures" (including the one in the next paragraph) involve  ``an incoherent group" $\RG$ satisfying the following condition (cf. \eqref{star0})
\begin{align}\label{star}
\text{\em $\RG(F\otimes_\BQ\BR)$ is compact modulo its center},
\end{align}
 a condition labeled as  
$(\ast)$ in Gan--Gross--Prasad \cite[\S27]{GGP}. Here, following Kudla's terminology, by ``an incoherent group" we refer to the isometry group associated to an incoherent collection of quadratic or Hermitian spaces.

More recently, there is another construction of special cycles based on a symmetric pair $(\RH,\RG)$, a (not necessarily inner) form of the pair $(\GL_n\times \GL_n,\GL_{2n})$, for example $(\RU(n)\times \RU(n),\RU(2n))$. See the introduction of \cite{LiQ}  for the relevant global Shimura varieties and subvarieties. There are new algebraic cycles on Shimura varieties associated to unitary groups which may locally be the unit group of any central simple algebra, a feature not seen in the AGGP conjecture or Kudla's program, see also the work in progress by Li--Mihatsch \cite{LM1,LM2}. For the $(\RU(n)\times \RU(n),\RU(2n))$-automorphic period integral (``unitary Friedberg--Jacquet period"), a relative trace formula approach is proposed and the stable FL on Lie algebra is established in \cite{LXZ}; the relative endoscopic FL on Lie algebra is also established in \cite{Les}. 


On the more speculative side, from \cite{Z23} we see that there may also be ``exceptional" examples, for which there is a local cycle (i.e., on a local Shimura variety) but there is no global Shimura variety.  Besides the general Bessel case $\U(m)\times \U(m+2r+1)$ and $\SO(m)\times \SO(m+2r+1)$, there is also the (local) cycle for the Ginzburg--Rallis $(H,G)$, which is related to the exterior cube L-function of automorphic representation on $\GL_6$. We will have to wait for the invention of new ``global Shimura varieties" to formulate the corresponding Gross--Zagier problem!

\subsection{Some arithmetic applications: the Beilinson--Bloch--Kato conjecture} 
\label{ss:BBK}
In the introduction, we have mentioned that Waldspurger's formula was a crucial ingredient in the proof of the  Birch--Swinnerton-Dyer conjecture for (modular) elliptic curves over a totally real number field in the analytic rank zero case, by Bertolini--Darmon's variant of the Kolyvagin system argument for Heegner points on a Shimura curve. In a higher dimensional (than curve) case, Yifeng Liu first proved the rank zero analog for the ``twisted triple product" automorphic motives in the product of a Shimura curve and a Hilbert--Blumenthal modular surface \cite{Liu16}, which may be placed in the framework of the triple product L-function and the Gross--Kudla--Schoen diagonal cycle (or equivalently  the $\SO(3)\times\SO(4)$-case of the Gross--Prasad construction); see also  \cite{Liu19,LT,W22} for works on other triple product L-functions. 

In \cite{LTXZZ} we extend this paradigm to the Rankin--Selberg motives in the context of  the Gan--Gross--Prasad construction for unitary groups $\RG=\RU(n-1)\times\RU(n)$. Let $\pi$ be a cohomological stably cuspidal tempered automorphic representation of $\RG(\BA)$ appearing in the cohomology \eqref{eq:decomp H}. Under certain technical conditions (mostly on the image of the residue Galois representation, see {\it loc. cit.}), we show that 
$$
L(1/2,\pi)\neq 0\imp  H^1_f(F', \rho_{\pi})=0.
$$
This confirms the corresponding Beilinson--Bloch--Kato conjecture in the analytic rank zero case. 

We also obtain a partial result towards the rank one case: under certain technical conditions on $\pi$, we show that {\em if the class  $\int_{z_{K_0}}\varphi\in H^1_f(F',\rho_\pi)$ in the Bloch--Kato Selmer group \eqref{eq: sel cl} is non-zero, then $\dim H^1_f(F', \rho_{\pi})=1$}. This is a direct generalization of Kolyvagin's theorem for Heegner points. Using a different Euler system argument (more close to Kolyvagin's Heegner point Euler system), 
Jetchev, Nekov\'a\v{r}, and Skinner also obtain a similar result (under different conditions) as well as one divisibility in a generalized Iwasawa Main conjecture.  
Combining the above results in the rank one case with the $p$-adic variant of AGGP in \S\ref{ss:pAGGP}, we would obtain a $p$-adic version of the Beilinson--Bloch--Kato conjecture formulated by Perrin-Riou  (under technical conditions on $\pi$)
$$
\ord L_{\it p}(1/2,\pi)=1\imp  \dim H^1_f(F', \rho_{\pi})=1.
$$
We refer to \cite{DZ} for more detail.

We would like to mention a series of recent remarkable works by Loeffler and Zerbes  \cite{LPSZ,LSZ,LZ}, partly joint with Pilloni and Skinner, in which they utilize the Gross--Prasad construction for orthogonal groups $\SO(4)\times \SO(5)$ and the exceptional isogenies $\SO(4)\sim \SL(2)\times \SL(2)$ and $\SO(5)\sim \Sp(4)$. As one of the applications, Loeffler and Zerbes \cite{LoZ} proved that, under certain technical hypotheses, for a modular abelian surface $A$ over $\BQ$,  
$$
\ord_{s=1}L(s,H^1(A))=0\imp \rank A(\BQ)=0,
$$
as predicted by (a generalization of the)  Birch--Swinnerton-Dyer conjecture! Rather than using Chow cycles as in the above work on Rankin--Selberg motives, their work makes use of ``higher Chow cycles" (elements in higher Chow groups), similar to Kato's Euler system argument using Siegel's units on modular curves (this seems to be  one of the reasons why the result is currently only known when the base field is $\BQ$). One of the difficulties is that the Tate module of an abelian surfaces has non-regular Hodge--Tate weights, hence does not appear directly in the (\'etale) cohomology of any Shimura varieties. To overcome this difficulty, a limit argument is necessary to extrapolate the Iwasawa theoretical Euler system from the regular case to the non-regular case.  

For another example of employing ``higher Chow cycles", Zhou \cite{Zhou} proved  a level-raising result using higher cycle classes on the special fiber of quaternionic Shimura varieties attached to an inner form of $\RR_{F/\BQ}\GL_2$ for a totally real field $F$. It would be interesting to know whether the local classes constructed in such a way can be lifted to global ones.

\subsection{Function field analogs }
\label{ss:fun fd}

When the global field $F$ is a function field,  there is an even further upgrade of the embedding $[H]\subset [G]$ to a morphism of moduli spaces of Shtukas \cite{Yun18,Yun22} with arbitrary number of legs (the Shimura variety case behaves similarly to the ``one-leg" case).  We refer to the article of Yun \cite{Yun24} for the definitions and the properties of the moduli spaces of Shtukas. In \cite{YZ1, YZ2}, in the case \eqref{eq:H G},
\begin{align*}
H=\RR_{F'/F}\BG_m, \quad G=B^\times,
\end{align*}Zhiwei Yun and the author express arbitrary higher order derivatives of L-functions of mildly ramified $\pi$ in terms of the intersection numbers of $\pi$-isotypical components of a special cycle associated to $H$. See also the related work by Howard and Shnidman \cite{HS,Sh}, and by C. Qiu \cite{Q}. 
In recent preprints in \cite{FYZ1,FYZ2,FYZ3}, Feng, Yun and the author investigated a function field analog of Kudla's program on generating series of special cycles and the arithmetic Siegel--Weil formula. We also refer to the article of Feng--Harris \cite[\S3,4]{FH} for some discussion related to the use of derived algebraic geometry in our works.

\bibliographystyle{amsalpha}

\end{document}